\documentclass[times,11pt]{article} 
\usepackage[small,compact]{titlesec}
\usepackage{verbatim}
\usepackage{wrapfig}
\usepackage{mathrsfs}
\usepackage{amsmath,amssymb,latexsym}
\usepackage{amsfonts}
\usepackage{times,lastpage,bm,xspace}
\usepackage{mathrsfs,bbm}
\usepackage{dsfont}
\usepackage{algorithm}
\usepackage{algorithmic}
\usepackage{stmaryrd}
\usepackage{booktabs}

\usepackage[pdftex]{graphicx}
\usepackage{url}

\usepackage[usenames]{color}
\definecolor{plum}  {rgb}{.6,0,.6}
\definecolor{forest}  {rgb}{0,.6,0}
\definecolor{midnight}  {rgb}{0,0,.8}
\usepackage[pdftex,plainpages=false,pdfpagelabels,colorlinks=true,urlcolor=midnight,linkcolor=midnight,citecolor=forest]{hyperref}

\usepackage[caption=false,font=footnotesize]{subfig}

\usepackage{amsmath, amsfonts, amssymb, amsthm,color,bm}

\def\ie{{\em i.e.,~}}
\def\eg{{\em e.g.,~}}
\def\cf{{\em cf.,~}}

\def\deq{\triangleq}
\DeclareMathOperator{\pen}{pen}
\DeclareMathOperator{\KL}{KL}
\DeclareMathOperator*{\argmin}{argmin}

\newcommand{\ave}[1]{\langle #1 \rangle}

\def\bI{\ones}
\def\cC{\mathcal{C}}
\def\cF{\mathcal{F}}

\def\cH{\mathcal{H}}
\def\cG{\mathcal{G}}

\newcommand{\wh}[1]{{\widehat{#1}}}
\newcommand{\wt}[1]{{\widetilde{#1}}}
\def \reals{\mathbb{R}}
\def \expect{\mathbb{E}}
\def \prob{\mathbb{P}}
\def\bI{\mathbbm{1}}

\newtheorem{lemma}{Lemma}[section]
\newtheorem{definition}{Definition}[section]
\newtheorem{theorem}{Theorem}[section]

\newtheorem{assumption}{Assumption}[section]
\newtheorem{remark}{Remark}[section]

\newcommand{\nm}{n}
\newcommand{\nc}{p}
\newcommand{\inty}{T}
\newcommand{\bD}{\bar{D}}
\newcommand{\btheta}{\bar{\theta}}
\newcommand{\FunClass}{\mathcal{F}_{\nc,s,D}}
\newcommand{\T}{\top}
\newcommand{\rip}{\delta_{s}}
\newcommand{\Ru}{(1 + \rip)}
\newcommand{\Rl}{(1 - \rip)}
\newcommand{\si}{\theta^i}
\newcommand{\sj}{\theta^j}

\def\eps{\frac{1}{2n}}
\def\shft{\frac{a_u-2a_\ell}{\sqrt{\nm}}} 
\def\resc{2{(a_u-a_\ell)\sqrt{\nm}}} 

\begin{document}
	\begin{center}
	{\bf{\LARGE{Minimax Optimal Rates for Poisson Inverse Problems with Physical Constraints}}}

	\vspace*{.1in}
	\begin{tabular}{ccc}
	Xin Jiang$^{1}$ & Garvesh Raskutti$^{2}$ & Rebecca Willett $^{3}$ \\
	\end{tabular}

	\vspace*{.1in}

	\begin{tabular}{c}
	  $^1$ Department of Electrical and Computer Engineering, Duke University\\
	  $^2$ Department of Statistics, University of Wisconsin-Madison\\
	  $^3$ Department of Electrical and Computer Engineering, University of Wisconsin-Madison
	\end{tabular}

	\vspace*{.1in}


	\end{center}

\begin{abstract}
This paper considers fundamental limits for solving sparse inverse
problems in the presence of Poisson noise with physical
constraints. Such problems arise in a variety of applications,
including photon-limited imaging systems based on compressed
sensing. Most prior theoretical results in compressed sensing and
related inverse problems apply to idealized settings where the noise
is i.i.d., and do not account for signal-dependent noise and physical
sensing constraints. Prior results on Poisson compressed sensing with
signal-dependent noise and physical constraints in~\cite{dense_pcs}
provided upper bounds on mean squared error performance for a specific
class of estimators. However, it was unknown whether those bounds were
tight or if other estimators could achieve significantly better
performance. This work provides minimax lower bounds on mean-squared error for
sparse Poisson inverse problems under physical constraints. Our lower
bounds are complemented by minimax upper bounds. Our upper and lower
bounds reveal that due to the interplay between the Poisson noise
model, the sparsity constraint and the physical constraints: (i) the
mean-squared error does not depend on the sample size $n$ other than
to ensure the sensing matrix satisfies RIP-like conditions and the
intensity $T$ of the input signal plays a critical role; and (ii) the
mean-squared error has two distinct regimes, a low-intensity and a
high-intensity regime and the transition point from the low-intensity
to high-intensity regime depends on the input signal $f^*$. In the
low-intensity regime the mean-squared error is independent of $T$
while in the high-intensity regime, the mean-squared error scales as
$\frac{s \log p}{T}$, where $s$ is the sparsity level, $p$ is the
number of pixels or parameters and $T$ is the signal intensity.  
\end{abstract}

\section{Introduction}

In this paper we investigate minimax rates associated with Poisson
inverse problems under physical constraints, with a particular focus
on compressed sensing (CS) \cite{CS:candes2,CS:donoho}. The idea
behind compressed sensing is that when the underlying signal is sparse
in some basis, the signal can be recovered with relatively few
sufficiently diverse projections. While current theoretical CS
performance bounds are very promising, they are often based on
assumptions such as additive signal-independent noise which do not
hold in many realistic Poisson inverse problem settings.  

Poisson inverse problems arise in a variety of applications.
\begin{itemize}
\item Several imaging systems based on compressed sensing have been proposed. 
The Rice single-pixel camera \cite{riceCamera} or structured illumination fluorescence microscopes \cite{CandesSIM} collect projections of a scene sequentially, and collecting a large number of observations within a limited time frame necessarily limits the number of photons that are detected 
 and hence the signal-to-noise ratio (SNR). In particular, if we collect $\nm$ projection measurements of a scene with $\nc$ pixels over a span of $\inty$ seconds, then the average number of photon counts per measurement scales like $\inty/\nm$. Thus even if the hardware used to determine which projections are measured (\eg a digital micromirror device) is capable of operating at very high speeds, we will still face photon limitations and hence Poisson noise \cite{snyder}.
\item In conventional fluorescence microscopy settings (\ie non-compressive, direct image acquisition) we may wish to represent a signal using a sparse superposition of Fourier and pixel elements to facility quantitative tissue analysis \cite{nimmi_spie,mca,gmca}. In this context one may wish to characterize the reliability of such a decomposition to ensure sufficient acquisition time or sufficient quantities of fluorescent protein. In this context $\nm$ is the number of pixels observed, and $\nc$ is the number of coefficients to be estimated.
\item In network flow analysis, we wish to reconstruct average  packet arrival rates  and instantaneous packet
counts for a given number of streams (or flows) at a router in a communication network,  where the arrivals of
packets  in each  flow are assumed to follow  a Poisson  process. All packet counting must be done in hardware at the router, and
any hardware implementation must strike a delicate balance between speed, accuracy, and cost \cite{ElephantsMice,CounterBraids,fishing,expander_pcs}. In this context, $\nm$ is the number of (expensive) memory banks in a router and ideally will be much less than $\nc$, the total number of flows to be estimated. $\inty$ corresponds to the total number of packets observed. 
\item Poisson inverse problems also arise in DNA analysis \cite{LaureDNA}, pediatric computer-aided tomography \cite{boone2003dose}, Twitter data analysis \cite{socioscope}, and quantitative criminology \cite{PoissonCrime}.
\end{itemize}
In all the above settings, when small numbers of events (\eg photons, x-rays, or packets) are observed, the noise or measurement error can be accurately  modeled using a Poisson distribution \cite{snyder}, so recovering a signal from projection measurements amounts to solving a Poisson inverse problem. Furthermore, the Poisson inverse problem is subject to physical constraints (e.g. hardware constraints, flow constraints) which are outlined below. In general, the available data acquisition time places fundamental limits on our ability to accurately reconstruct a signal in a Poisson noise setting. Until now, relatively little about the nature of these limits was understood.

\subsection{Observation model and physcial constraints}

In this paper we consider an observation model of the form
\begin{equation}
\begin{aligned}
y\sim {\rm Poisson}(\inty Af^*),
\end{aligned}
\label{eq:model}
\end{equation}
where $A \in \reals_+^{\nm \times \nc}$ in a sensing matrix corresponding to the $\nm$ different projections of our signal of interest $f^* \in \reals^\nc_+$ and $\inty \in \reals_+$ is the total data acquisition time. In particular, \eqref{eq:model} is a shorthand expression for the model 
\begin{equation}
y_i {\sim} {\rm Poisson}\left(\inty \sum_{j=1}^\nc A_{i,j}f^*_j\right), \qquad i \in \{1,\ldots,\nm\},
\notag
\end{equation}
where the $y_i$'s are independent.
Since for the Poisson distribution the variance depends on the signal intensity, the noise of the observed signal $y$ depends on $T$, $A$ and $f^*$. 

In particular, $f^*$ corresponds to the rate at which events are being generated, which are intrinsically nonnegative. In addition,  $A$ 
must be composed of nonnegative real numbers with each column summing to at most one. Specifically, $A$ must
satisfy the following physical constraints:
\begin{subequations}
\begin{align}
A &\succeq 0 && \mbox{nonnegative}; \label{eq:LS1}\\
\bI_{\nm \times 1}^\T A & \preceq \bI_{\nc \times 1} && \mbox{bounded column sums}.\label{eq:LS3}
\end{align}
\end{subequations}
Together these constraints imply (a) $A_{i,j}  \in [0,1] \; \forall (i,j)$ and (b) $\|Af\|_1 \le \|f\|_1$ for all nonnegative $f$; this latter property with equality is often called {\em flux preservation}, particularly in the astronomy community. Thoughtout this paper, we also refer to constraint~\eqref{eq:LS3} as the flux-preserving constraint.

To build intuition about these constraints it is helpful to consider an imaging system in which $A_{i,j}$  corresponds to the likelihood of a photon generated at location $j$ at the source $f^*$ hitting our detector at location $i$. Such a likelihood is necessarily nonnegative and bounded, and the $j^{\rm th}$ column sum, which corresponds to the likelihood of a photon from location $j$ hitting any of the $\nm$ detectors, must be bounded by one. 
Without the constraint \eqref{eq:LS3}, we could have an unrealistic imaging system in which more photons are detected than we have incident upon the aperture. These constraints arise naturally in all the motivating applications described above using similar arguments. 

We assume $\|f^*\|_1 = 1$ throughout, so the total expected number of observed events is proportional to $\inty$; thus $\inty$ alone (\ie neither $f^*$ nor $A$) controls the signal-to-noise ratio. This model is equivalent to assuming the signal of interest integrates to one (like a probability density), and by allotting ourselves $\inty$ units of time to acquire data (\eg collect photons). 

To incorporate sparsity, we assume $f^*$ is $(s+1)$-\emph{sparse} in a basis spanned by the columns of the orthonormal matrix $D \in \reals^{\nc \times \nc}$, where we write $D = [d_1, \ldots, d_{\nc}]$ for $d_i \in \reals^\nc$ for all $i$. We assume that $d_1 = \nc^{-1/2}\bI_{\nc \times 1}$. \footnote{This assumption precludes signals that are sparse in the canonical or pixel basis. Our proof techniques, with some technical adjustments, can be adapted to that setting.}
Standard choice of basis matrices $D$ are Fourier and wavelet basis matrices which we discuss in more detail later. 
For an arbitrary vector $v \in \reals^{\nc}$, let $\|v\|_0 \deq \sum_{j=1}^{\nc}{\mathbf{1}(v_j \neq 0)}$ denote the standard $\ell_0$-norm.

Hence, if we define $\theta^* \deq D^\T f^* \in \reals^\nc$, then $\|\theta^*\|_0 \le s+1$. Note that by construction, $\theta^*_1 \equiv \ave{d_1,f^*} \equiv 1/\sqrt{\nc}$. Let $\bD \deq [d_2, d_3, \ldots, d_\nc] \in \reals^{\nc \times (\nc-1)}$ be all but the first (constant) basis vector,  let $\btheta^* \deq \bD^\T f^*$ be all but the first (known) coefficients, and note that $\btheta^*$ is $s$-sparse.

To summarize, we consider $f^*$ belonging to the set: 
\begin{equation}
\label{eqn:funClass}
\FunClass \deq \left\{f \in \reals_+^{\nc} \; : \;  \|f\|_1  = 1,~\|\bD^\T  f\|_0 \leq s \right\}.
\notag
\end{equation}
We assume throughout the paper that $f^* \in \FunClass$. 

\subsection{Our Contributions}

In this paper, we derive mean-squared error bounds for the Poisson inverse problem described in ~\eqref{eq:model} where the function $f^*$ belongs to $\FunClass$. Our main results are presented in Section~\ref{sec:Results}. In particular, in Section~\ref{sec:LowerBound} we provide a minimax lower bound for the mean-squared error, restricted to the set of signals and sensing matrices that satisfy the above physical constraints. In Section~\ref{sec:UpperBound} we provide a minimax upper bound based on analysis of a modified version of the estimator considered in~\cite{dense_pcs}.

Our upper and lower bounds reveal two new effects due to the interplay between the Poisson model, the sparsity assumption and the physical constraints. The first new effect is that the flux-preserving constraint ensures that the mean-squared error has no dependence on the sample size $n$ outside of needing $n$ to be sufficiently large for $A$ to satisfy RIP-like assumptions. Instead the intensity $T$ plays a significant role in the mean-squared error. Secondly, the interplay between the non-negativity constraint and the sparsity assumption means that there are two distinct regimes, a \emph{low-intensity} and \emph{high-intensity} regime. In the low-intensity regime, \ie for $\inty$ below some critical level, the mean-squared error is independent of $T$, whereas in the high-intensity regime the mean-squared error scales as $\frac{s \log p}{T}$, which we prove is the minimax optimal rate. The point at which the mean-squared error transitions from the low-intensity to the high-intensity regime depends significantly on the input signal $f^*$ and the orthogonal basis matrix that induces sparsity on $f^*$. Hence our upper and lower bounds are sharp up to constant in the high-intensity regime whereas in the low-intensity regime, our upper and lower bounds do not match since the optimal choices of $f^*$ depend significantly on the orthonormal basis matrix $D$. The theoretical performance bounds are supported by a suite of simulations presented in Section~\ref{sec:Discussion}. Thus the main theoretical results in this paper are perhaps surprising and certainly at odds with similar results appearing in Gaussian or bounded noise settings where the bounds are independent of the sparsifying basis and the choice of $f^*$, and the mean-squared error bounds depend on $n$.

\section{Assumptions and Main Results}
\label{sec:Results}

We provide upper and lower bounds for the minimax mean-squared error. Let $\wh{f} \equiv \wh{f}(y)$ denote an estimator of $f^*$ (\ie $\wh{f}$ is a measurable function of $y$, and $\inty AD$). The performance of the estimator is evaluated using the following  squared $\ell_2$ risk:
\begin{equation}
\mathcal{R}(\wh{f}, f^*) \deq \expect[\| \wh{f}-f^*\|_2^2],
\notag
\end{equation}
where the expectation is taken over $y$. The minimax risk is defined as:
\begin{equation}
\min_{{\wh{f}}}\max_{f^* \in \FunClass} \mathcal{R}(\wh{f}, f^*).
\notag
\end{equation}
Here the $\min$ is taken over the set of measurable functions of $(\inty AD,y)$. To begin we state assumptions imposed on the sensing matrix $A$.

\subsection{Assumptions on $A$}
\label{sec:Assumptions}

\begin{assumption}\label{as:wtA}
There exist constants $a_\ell$ and $a_u$ such that $a_\ell < a_u$ and a matrix $\wt{A} \in [a_\ell/\sqrt{\nm},a_u/\sqrt{\nm}]^{\nm\times\nc}$
\begin{equation}
A=\frac{\wt{A}+\shft\bI_{\nm\times\nc}}{\resc}.
\label{eq:Aconstruct}
\end{equation}
\end{assumption}

\begin{assumption}\label{as:Arip1}
For all $u\in \reals^{\nc}, \|u\|_0\le 2s$, 
there exists a constant $\rip \equiv \rip(\nm,\nc) > 0$ such that
\begin{equation}
\|\wt{A}D u\|_2^2 \le (1+\rip) \|u\|_2^2.
\notag
\end{equation}
\end{assumption}

\begin{assumption}\label{as:Arip2}
For all $u\in \reals^{\nc}, \|u\|_0\le 2s$, there exists a constant $\rip\equiv \rip(\nm,\nc)>0$ such that
\begin{equation}
\|\wt{A}D u\|_2^2 \ge (1-\rip) \|u\|_2^2.
\notag
\end{equation}
\end{assumption}
Assumption~\ref{as:wtA} ensures that $A \succeq 0$. Note that in many practical cases in which Assumptions~\ref{as:wtA},~\ref{as:Arip1}, and~\ref{as:Arip2} are satisfied, we will have  $a_\ell < 0 < a_u$.
Assumptions~\ref{as:Arip1} and~\ref{as:Arip2} together are the the RIP conditions for $\wt{A}D$ introduced in~\cite{CandesTao05}. 
Our minimax lower bound requires only Assumption~\ref{as:Arip1} whereas our upper bound requires both Assumption~\ref{as:Arip1} and Assumption~\ref{as:Arip2}, which is consistent with~\cite{garvesh}. 

Assumption~\ref{as:wtA} also introduces a new matrix $\wt{A}$. This reflects the actual process of constructing a sensing matrix that satisfies the Conditions~\eqref{eq:LS1} and~\eqref{eq:LS3}. The typical design process would be to start with a bounded matrix which satisfies the RIP conditions, rescale each entry and add an offset to the original matrix to make sure the new matrix is positive and flux-preserving. The following lemma ensures that these two constraints are satisfied via such operations:

\begin{lemma}
\label{lem:flux}
If the matrix $\wt{A}$ satisfies $\wt{A}_{i,j}\in\frac{1}{\sqrt{\nm}}[a_\ell,a_u]$ for all entries, and the sensing matrix $A$ is generated according to Eq.~\eqref{eq:Aconstruct}, then $A$ satisfies Conditions~\eqref{eq:LS1} and~\eqref{eq:LS3}.
\end{lemma}
\noindent The proof is provided in Section~\ref{sec:proofflux}.

For example, let $A$ be constructed based on a shifted and rescaled Bernoulli ensemble matrix $\wt{A}$, where
 \begin{equation}
\begin{aligned}
\prob(\wt{A}_{i,j})=
\begin{cases}
	1/2,& \wt{A}_{i,j} = -\sqrt{\frac{1}{\nm}} \\
	1/2,& \wt{A}_{i,j} = \sqrt{\frac{1}{\nm}}.
\end{cases}
\end{aligned}
\label{eq:A0construct2}
\end{equation}
The resulting sensing matrix $A$ has i.i.d. entries valued $1/(2\nm)$ or $1/\nm$ and using results in ~\cite{Baraniuk08},  $\wt{A}$ and $\wt{A}D$ satisfy the RIP conditions in Assumptions \ref{as:Arip1} and \ref{as:Arip2} for $\nm\ge C_0 s\log\left(\nc/s\right)$ with probability at least $1-e^{-C_1\nm}$.

\subsection{Lower bounds on minimax risk}
\label{sec:LowerBound}

In this section we present a lower bound on the minimax risk. The interaction between the orthonormal basis matrix $D$ and the sparsity constraint has an effect on the lower bound which is captured by the following $s$-sparse localization quantity:

\begin{definition}[$s$-sparse localization]  \label{def:localization} $\lambda_s$ is said to be the $s$-sparse localization quantity of a matrix $X$ if
\begin{equation}
\lambda_s = \lambda_s(X) \deq \max_{\substack{v \in \{-1,0,1\}^p \\ \|v\|_0 = s}} \|Xv\|_\infty.
\notag
\end{equation}
\end{definition}
The name ``localization'' in Definition~\ref{def:localization} derives
from a similar quantity defined in~\cite{patricia}. The key difference
between our localization constant and the one in \cite{patricia} is
that we incorporate the sparsity level directly within the
definition. Our minimax lower bound depends on $\lambda_{k}(\bD)$ for $1 \leq k \leq s$.

\begin{theorem}[Minimax lower bound]
\label{thm:lower}
For $k = 1,\ldots,s$, let $\lambda_{k} = \lambda_{k}(\bD)$ be the $k$-sparse localization quantity of $\bD$.
 If $\nc\ge10$, $1\le s < p/3-1$, and Assumptions~\ref{as:wtA} and~\ref{as:Arip1} hold with $0\le \rip < 1$,  
then there exists a constant $C_L>0$ that depends only on $a_u$ and $a_\ell$ such that 
\begin{equation}
\begin{aligned}
\min_{{\wh{f}}}\max_{f^*\in\FunClass} \mathcal{R} (\wh{f}, f^*)
\ge C_L\max_{1\le k\le s} \left\{ \min \left( \frac{k}{\nc^2\lambda_k^2}, \frac{k}{\Ru\inty}\log{\frac{\nc-k-1}{k/2}} \right) \right\}.
\end{aligned}
\label{eq:LB}
\end{equation}
\end{theorem}
\noindent The proof is provided in Section~\ref{sec:proofLower} and involves adapting the information-theoretic techniques developed in~\cite{garvesh,Han94,IbrHas81} to the Poisson noise setting. 

\begin{remark}
Within the the $\max$ term, the lower bound is a minimum of two terms. 
These two terms represent the minimax rate in a ``low-intensity'' regime where $\inty$ is less than $\displaystyle O({ \min_{k\le s} \lambda_{k}^2}\nc^2\log \nc )$ and a ``high-intensity" regime where $\inty$ is above $\displaystyle O({ \max_{k\le s} \lambda_{k}^2} \nc^2 \log \nc )$. 
The sparse localization quantity $\lambda_k$ plays a significant role in the performance at lower intensities, whereas at high intensities the scaling is simply $O(\frac{s \log\nc}{\inty})$. 
\end{remark}

\begin{remark}
It should be noted that the quantity $\lambda_k$ is not a constant in general. It is a function of $D$ and $k$, which also means that $\lambda_k$ can be a function of the signal dimension $\nc$. Depending on the choice of the basis,  $\lambda_k$ can be on the order of $k/\sqrt{\nc}$ (\eg for  the DCT basis), or a constant (\eg for the Haar wavelet basis). The computation of $\lambda_k$ for various $D$  is  described in Section~\ref{sec:UpperLower}. As a  result, our lower bound in the low-intensity setting for the discrete cosine basis is different from the Haar wavelet basis (see Table~\ref{tab:lowint}). 
\end{remark}

\begin{remark}
The bound in \eqref{eq:LB} is independent of $n$, but note that the assumptions only hold for $n$ sufficiently large. For instance, if $A$ were generated using a Bernoulli ensemble as in \eqref{eq:A0construct2}, we would need $\nm = O(s\log(\nc/s))$ to ensure Assumptions~\ref{as:Arip1} and~\ref{as:Arip2} hold. 
\end{remark}

\subsection{Upper bound}
\label{sec:UpperBound}
Now we provide an upper bound on the minimax rate. To obtain the upper bound, we analyze an $\ell_0$-penalized likelihood estimator $\wh{f}$ introduced in ~\cite{dense_pcs}.

\begin{theorem}[Minimax upper bound]
\label{thm:upper}
If $f^* \in \FunClass$, $A$ satisfies Assumption~\ref{as:wtA},~\ref{as:Arip1}, and~\ref{as:Arip2} with $0\le \rip< 1$, then there exists a constant $C_U>0$ depending only on $a_u$ and $a_\ell$ such that 
\begin{equation}
\begin{aligned}
\min_{\wh{f}} \max_{f^* \in \FunClass} \mathcal{R} (\wh{f}, f^*)
	 	 \le C_U \min \biggr(
	 \frac{{s}\log_2\nc}{\Rl\inty} 
	 + \frac{{s}\log_2(\inty+1)}{\Rl\inty}, s \max_{j,k}|D_{j,k}|^2, 1\biggr). 
\end{aligned}
\notag
\end{equation}
Additionally, as long as $\inty$ does not increase exponentially with $\nc$, ${\log_2(\inty+1)}$ is dominated by $c \log_2\nc$ for some $c > 0$, our upper bound satisfies
\begin{equation}
\begin{aligned}
\min_{\wh{f}} \max_{f^* \in \FunClass} \mathcal{R} (\wh{f}, f^*)
	 	 \le C'_U \min \biggr(
	 \frac{{s}\log_2\nc}{\Rl\inty}, s \max_{j,k}|D_{j,k}|^2, 1 \biggr).
\end{aligned}
\label{eq:upper}
\end{equation}
\end{theorem}
The proof is provided in Section~\ref{sec:proofUpper}.

\begin{remark}
The first term $\frac{s \log p}{T}$ arises from analyzing a penalized maximum likelihood estimator. The second and third term arises from using the $0$ estimator and upper bounding $\| f^*\|_2^2$ by $\min( s \max_{j,k}|D_{j,k}|^2, 1)$.
\end{remark}

\begin{remark}
Like the lower bound, we see different behavior of the reconstruction
error depending on the intensity $T$. For $T > \min(s \log p,
(\max_{i,j}|D_{j,k}|)^{-2} \log p )$ the upper bound scales as
$\frac{s \log p}{T}$ while for $T \leq \min(s \log p,
(\max_{j,k}|D_{j,k}|)^{-2} \log p )$ the upper bound is independent of
$T$. Note that in the high-intensity regime (large $T$), the upper and
lower bounds of $\frac{s \log p}{T}$ match up to a constant. We
discuss the scaling in the low-intensity regime in more detail in
Section~\ref{sec:UpperLower}. 
\end{remark}

\begin{remark}
The lower bound in Theorem~\ref{thm:lower} only requires Assumptions~\ref{as:wtA} and~\ref{as:Arip1}, whereas the upper bound requires the additional Assumption~\ref{as:Arip2}. This is consistent with assumptions required for the minimax upper and lower bounds in ~\cite{garvesh}. 
\end{remark}

\begin{remark}
Like in the lower bound, the number of observations $\nm$ does not appear in the mean-squared error expression. 
\end{remark}

\section{Consequences and discussion}\label{sec:Discussion}

The minimax lower bounds show that the mean-squared error scales as $\frac{s \log(\nc/s )}{\inty}$ provided $\inty$ is sufficiently high and the upper bound scales as $\frac{s \log \nc}{\inty}$, where $\inty$ corresponds to the data acquisition time or the total number of event available to sense. 
Notice that the number of observations $\nm$ does not appear in the mean-squared error expression. Hence increasing the number of observations $\nm$ beyond what is necessary to satisfy the necessary assumptions on $A$ is not going to reduce the mean-squared error.
Furthermore, our analysis proves that the mean-squared error is inversely proportional to the intensity $\inty$, so $\inty$ plays the role that $\nm$ usually does in the standard compressed sensing setup.
We also carefully characterize the mean-squared error behavior when the intensity $\inty$ is low, and show {\em bounds that depend on the interactions between the sparsifying basis and the physical constraints.} We elaborate on these concepts below.

\subsection{High and low intensities}\label{sec:UpperLower}

As mentioned earlier, the lower bounds in Theorem~\ref{thm:lower} and the upper bounds in Theorem~\ref{thm:upper} agree up to a constant when the intensity is sufficiently high. In this section, we explore the relationship between the upper and lower bounds in the low-intensity setting for discrete cosine transform (DCT), discrete Hadamard transform (DHT) and discrete Haar wavelet basis (DWT). The elements for the DCT basis matrix are:

\begin{equation}
\begin{aligned}
D^{\rm DCT}_{j,k}=\begin{cases}
\sqrt{\frac{1}{\nc}}, & k=1, j=1,\ldots,\nc \\
\sqrt{\frac{2}{\nc}}\cos{\left(\frac{(2j-1)(k-1)\pi}{2\nc}\right)}, & k=2,\ldots,\nc, j=1,\ldots,\nc.
\end{cases}
\end{aligned}
\label{eq:dctval}
\end{equation}

The DHT matrix $D^{\rm DHT}_{(m)}\in\reals^{2^m\times 2^m}$ has the form
\begin{equation}
\begin{aligned}
\left(D^{\rm DHT}_{(m)}\right)_{j,k}=\frac{1}{2^{\frac{m}{2}}}(-1)^{j\cdot k}.
\end{aligned}
\label{eq:dhtval}
\end{equation}

The DWT matrix of dimension $2^m$ by $2^m$, the non-zeros have magnitudes in the set $\{2^{-m/2}, 2^{-(m-1)/2}, \ldots, 2^{-1/2}\}$. 

In the low-intensity setting, the lower bound scales like ${\displaystyle \max_{1\le k\le s}} \frac{k}{\nc^2\lambda_k^2}$, while the upper bound scales as $\min\left( 1, s\max_{j,k}|D_{j,k}|^2 \right)$. Since the bounds depend on $\lambda_k$ and $\max_{j,k}|D_{j,k}|$, it is not straightforward to see how these two bound match with each other. In Table~\ref{tab:lowint}, we show the comparison of the bounds for DCT basis, DHT basis and DWT basis. The calculation of $\lambda_k$ can be found in Section~\ref{sec:calcLambda}. 

\begin{table}[htbp]
\centering
\begin{tabular}{@{} ccccc @{}}
\toprule
		 		& $\lambda_k$
				& Lower bound		
				& $\max_{j,k} |D_{j,k}|$  
				& Upper bound \\
\midrule  
DCT	\& DHT		& $\frac{\sqrt{2}k}{\sqrt{\nc}}$	
				& $\min\left(\frac{1}{\nc},\frac{s\log{({\nc}/{s})}}{\inty}\right)$
				& $\frac{2}{\sqrt{\nc}}$ 
				& $\min\left(\frac{s}{\nc},\frac{s\log\nc}{\inty} \right)$\\
\midrule
DWT			 	& $\frac{1}{\sqrt{2}-1}$		
				& $\min\left(\frac{s}{\nc^2},\frac{s\log{({\nc}/{s})}}{\inty}\right)$
				& $\frac{1}{\sqrt{2}}$
				& $\min\left(1,\frac{s\log\nc}{\inty} \right)$\\
\bottomrule
\end{tabular}
\caption{Lower and upper bounds for three different bases (neglecting constants).}
\label{tab:lowint}
\end{table}

For the high-intensity setting, the lower and upper bound have the
same rate regardless of the basis (except for the minor difference in
the $\log$ term), which proves the bound is tight. For the
low-intensity setting, the upper and lower bounds do not match. 

To see why this is so, note that different sparse supports for
$\theta$ yield very different MSE performances because of the
interaction between the sparse support and the nonnegativity
constraints. Our proof shows that in low-intensity settings, the MSE is
proportional to the squared $\ell_2$ norm of the zero-mean signal 
$\|f^* - \bI_{n\times 1}/\sqrt{p}\|_2^2 = \|\btheta^*\|_2^2 $. Upper
and lower bounds on this quantity
depend on the interaction between the constraint, $\|f^*\|_1 = 1$, the sparsity constraint,
and the orthonormal basis matrix $D$. The general lower bound $\frac{k}{p^2 \lambda_k^2}$ derives from the following $\theta^*$:
 \begin{equation}
\begin{aligned}
\theta^*_j=
\begin{cases}
	\frac{1}{\sqrt{p}},& j=1 \\
	\pm\frac{1}{p \lambda_k},& 2 \leq j \leq k + 1\\
	0,& \text{otherwise}.
\end{cases}
\end{aligned}
\label{EqnThetaSim}
\end{equation}
It is straightforward to see that $\|\btheta^*\|_2^2 = 
\frac{k}{p^2 \lambda_k^2}$ and Section~\ref{sec:proofLower} shows that
$f^* = D \theta^*$ satisfies all the necessary constraints. On the
other hand, the upper bound $\max(1, s\max_{j,k}|D_{j,k}|^2 )$ follows from
\begin{equation}
\|\btheta^*\|_2^2 < \|\theta^*\|_2^2 = \|f^*\|_2^2 \leq \|f^*\|_1^2 = 1
\notag
\end{equation} 
and
\begin{equation}
\|\btheta^*\|_2^2 \leq s \|\theta^*\|_{\infty}^2 = s \|D^T f^*\|_{\infty}^2 \leq s \max_{j,k}|D_{j,k}|^2.
\notag
\end{equation}
Whether the upper or lower bounds are tight depend on the orthonormal basis matrix $D$. For example for the DWT basis, when $s$ is sufficiently large (\ie $\Omega(\log p)$), the
$\ell_2$ norm of $f \in \FunClass$ can be large, leading to tighter
lower bounds.
\begin{lemma}
\label{lem:tightLB}
Let $D$ be the discrete Haar wavelet basis matrix. If $D^\T f^*$ is $(s+1)$-sparse with $\nc\ge4$, $s\ge\log_2\nc$, $T\le \frac{(a_u-a_{\ell})^2\log\nc}{4\Ru}$, and Assumptions~\ref{as:wtA} and~\ref{as:Arip1} hold with $0\le \rip < 1$,  
then there exists an absolute constant $C_L''\ge1/8$ such that
\begin{equation}
\begin{aligned}
\min_{{\wh{f}}}\max_{f^*\in\FunClass}\expect\left[\|f^*-\wh{f}\|_2^2\right]
\ge C_L''.
\end{aligned}
\label{eq:tightLB}
\end{equation}
\end{lemma}
\noindent The proof is provided in Appendix~\ref{sec:prooftightLB}.

To provide an explanation for the effect of the matrix $D$ on the
minimax rates, the nonnegativity constraints implicitly impose limits
on the amplitudes of the coefficients of $f^*$, and these limits
strongly impact MSE performance. For instance, for 1-sparse $\btheta$
in the DWT basis, fine-scale basis vectors have amplitude of $O(1)$,
and so the corresponding basis coefficient cannot be too large
without violating the nonnegativity constraint; in contrast,
coarse-scale basis vectors have amplitude of $O(1/\sqrt{p})$ and can
have much larger coefficients without violating the nonnegativity
constraint. For the DWT basis, a sparse support can be chosen to
generate a high-$\ell_2$ norm signal that approximates a delta
function in $\FunClass$; this kind of signal dominates our upper
bounds. A different sparse support set may correspond to basis vectors
with mostly disjoint supports; such a signal must have a small
$\ell_2$ norm in order to be in $\FunClass$ (i.e. to be nonnegative);
this kind of signal dominates our lower bounds. A similar effect
occurs in the DCT basis.  We may contrast these behaviors with what
happens in a conventional Gaussian noise analysis without a
nonnegativity constraint, where we would have no constraints on
coefficient magnitudes to ensure nonnegativity.  This example and a
similar example for the DCT basis are presented with a simple
simulation described in Section~\ref{sec:supportSim}.

\subsection{Simulation results}
\label{sec:Simulations}

We further assess the proposed bounds by a series of experiments that examine the reconstruction performance of signals that satisfy our assumptions. We assume the signal is $s$-sparse under the DCT or DWT basis. Plots are not generated for the DHT, as it exhibits similar performance behavior to the DCT. 

We construct the signal $\theta^*$ as follows: the first coefficient
$\theta^*_1=1/\sqrt{p}$ corresponds to the DC level and ensures
$\|D\theta^*\|_1 = 1$. The locations of the remaining $s$ non-zero
coefficients are randomly generated, while the amplitude of each non-zero
coefficient is the same and chosen to be $\pm 1/(p\lambda_s)$, which
ensures that
$D\theta^*\in\FunClass$. (In particular, this construction of $\theta^*$
corresponds to elements of the packing set described in Section~\ref{sec:proofLower}.)
The sensing matrix $A$ is
generated according to Eqs.~\eqref{eq:Aconstruct}
and~\eqref{eq:A0construct2}. The reconstruction is done by solving a
$\ell_1$-penalized Poisson likelihood optimization problem. The reason for using a relaxed $\ell_1$ penalty instead of the $\ell_0$ penalty which is used for the analysis of the upper bound is that the $\ell_0$ penalization function is non-convex, and cannot be solved
efficiently. The optimization problem is solved using the SPIRAL
algorithm developed in \cite{zac}. Each data point in the plot shows
the  mean-squared error (MSE) averaged over 100 experiments.  

\subsubsection{Impact of sparsifying dictionary on rates}
\label{sec:supportSim}

As mentioned above, the the nonnegativity constraints
impose limits on the amplitudes of the coefficients of $f^*$, and
these limits strongly impact MSE performance. We illustrate that
effect here. In particular, we see that randomly selected sparse
supports and uniform coefficient magnitudes result in MSE rates which
correspond to our lower bounds. In contrast, carefully constructed
sparse supports (specific to the underlying sparsifying dictionary)
and highly non-uniform coefficient magnitudes result in slower MSE
rates which correspond to our upper bounds.

\begin{figure}[h!]
\centering
\subfloat[Sparse DCT coefficients]
  {\includegraphics[width=0.4\textwidth]{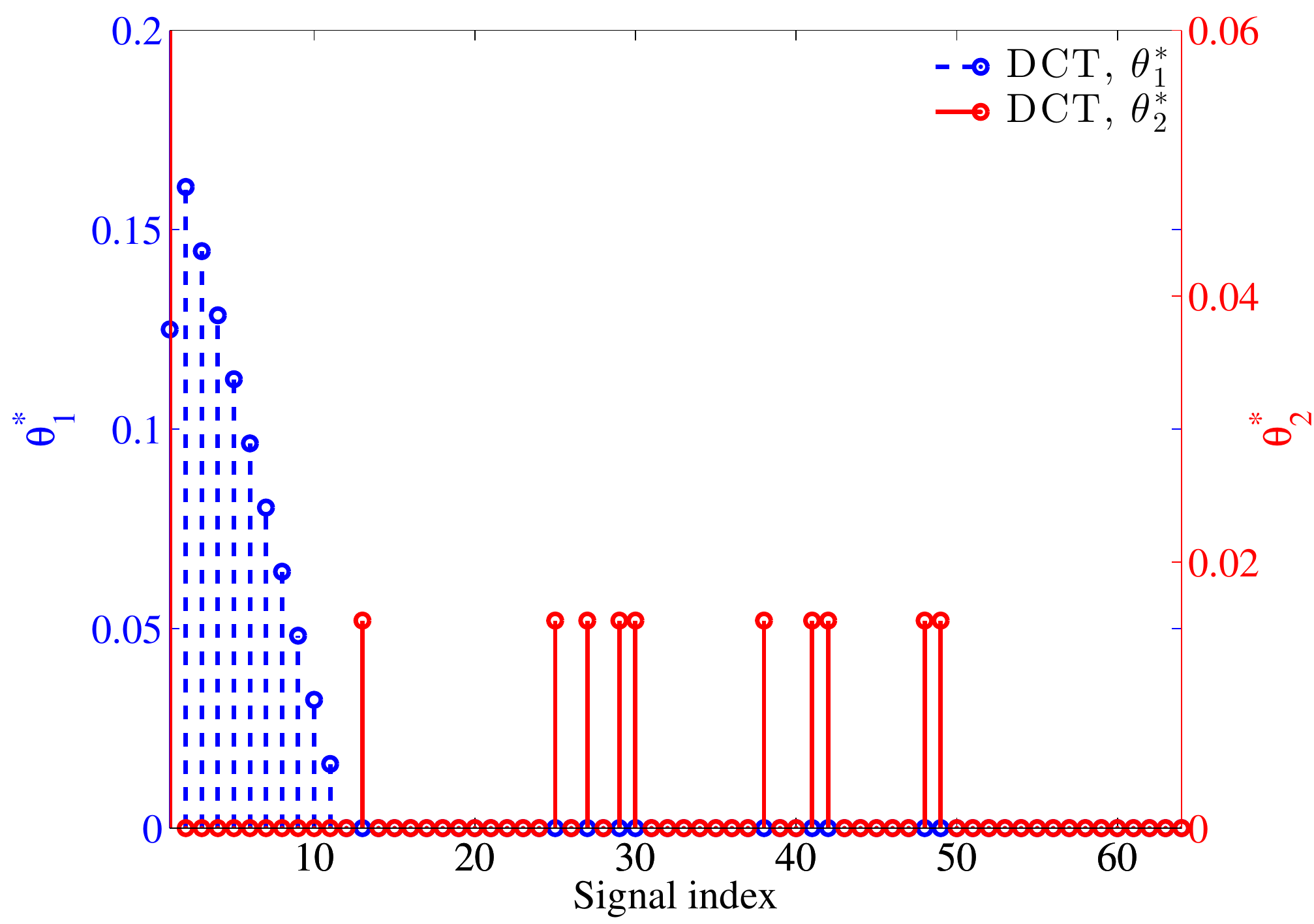}
  \label{fig:DCT_theta}}~
\subfloat[Sparse DWT coefficients]
  {\includegraphics[width=0.4\textwidth]{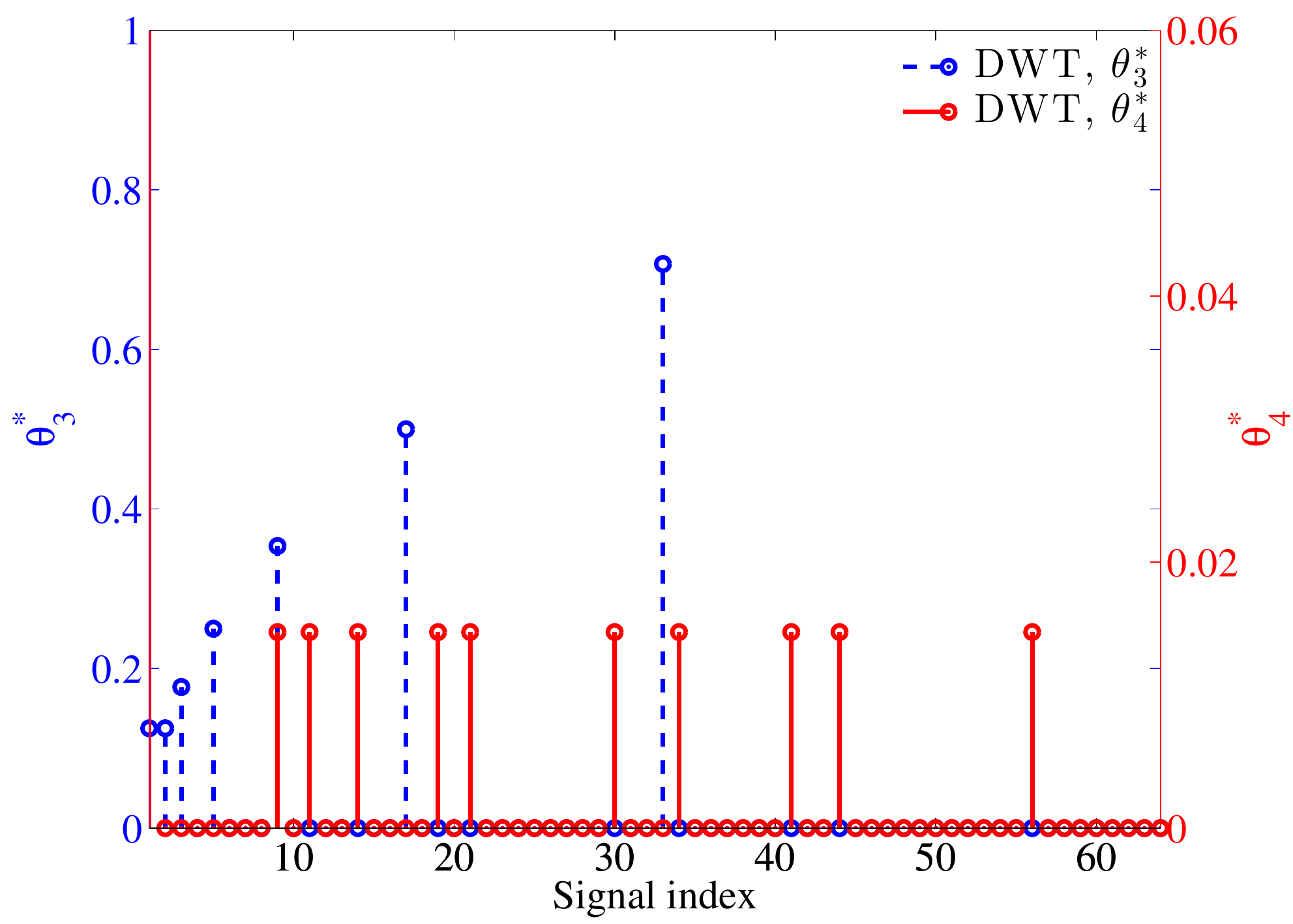} 
  \label{fig:DWT_theta}}\\
\subfloat[$f^*$ sparse in DCT Basis]
  {\includegraphics[width=0.4\textwidth]{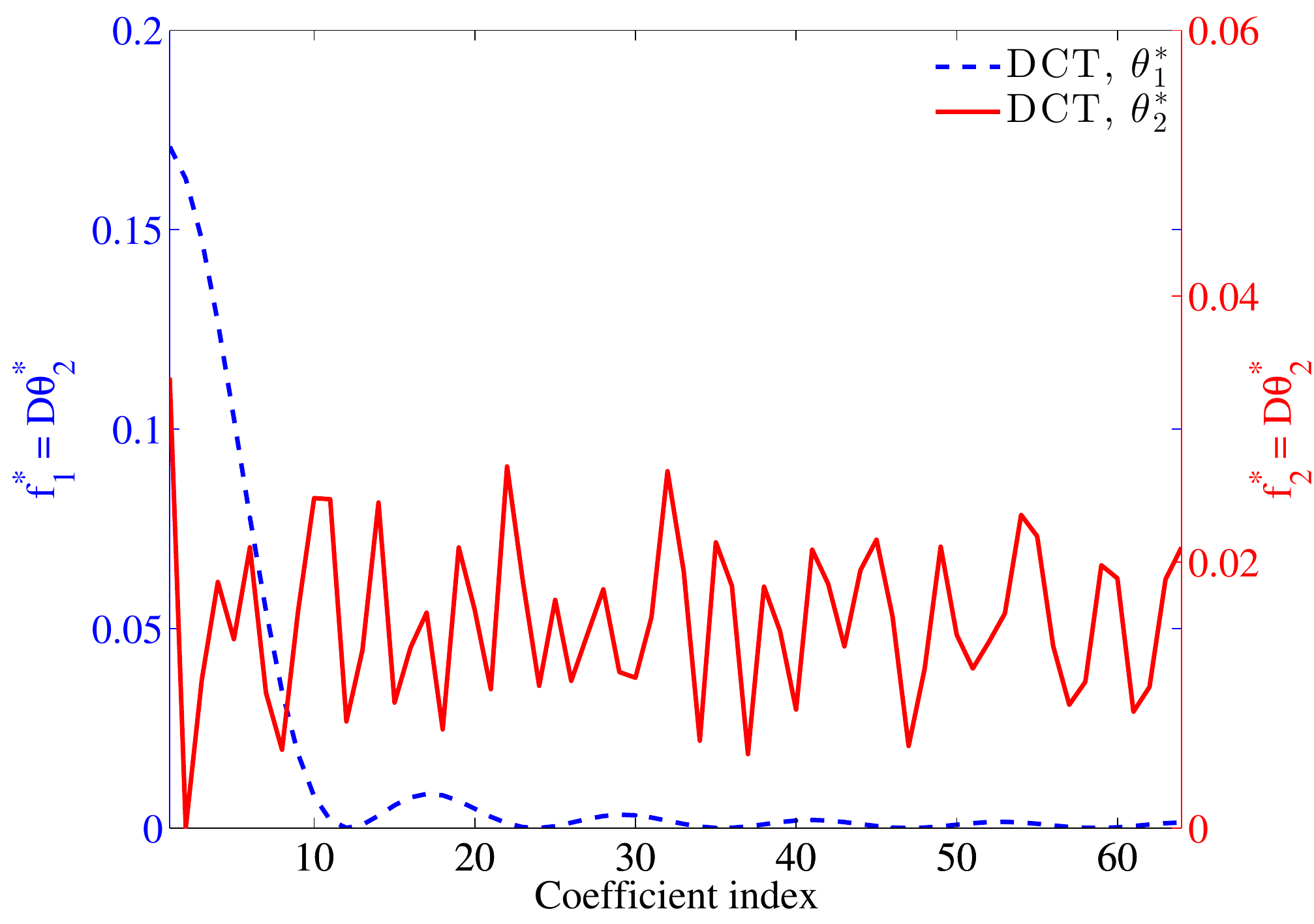}
  \label{fig:DCT_f}}~
\subfloat[$f^*$ sparse in DWT Basis]
  {\includegraphics[width=0.4\textwidth]{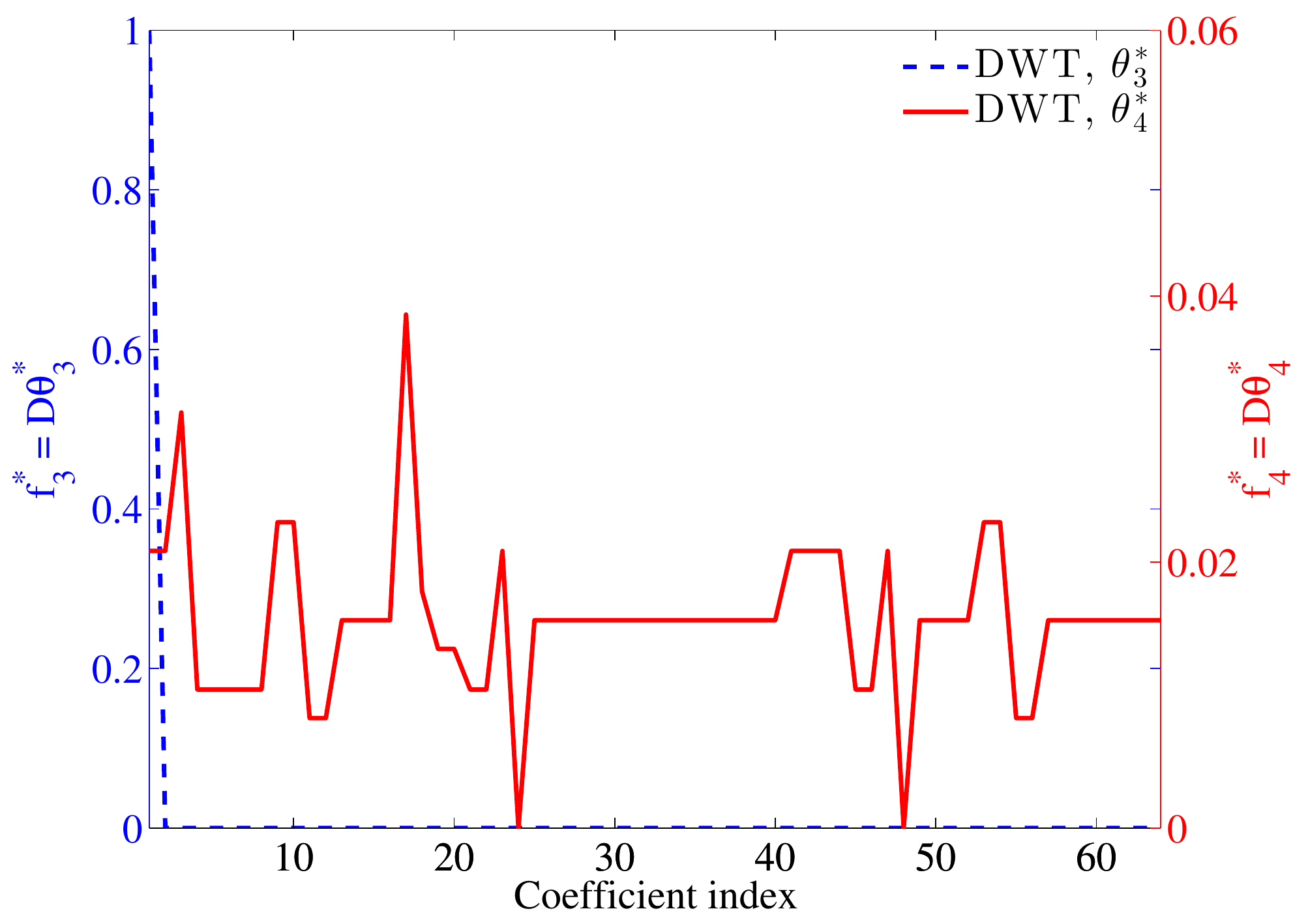} 
  \label{fig:DWT_f}}\\
\subfloat[MSE of DCT]
  {\includegraphics[width=0.4\textwidth]{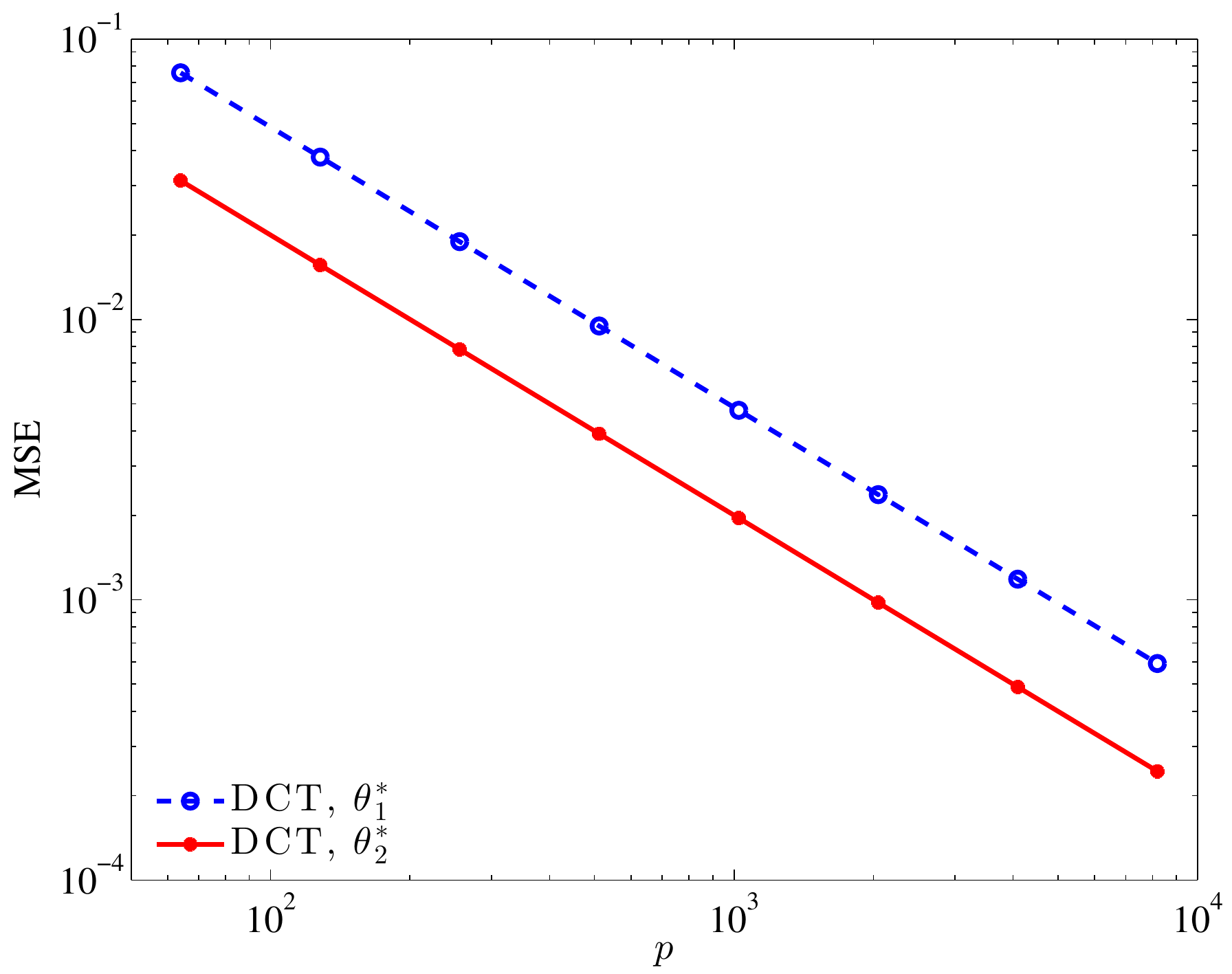}
  \label{fig:SP_DCT}}~
\subfloat[MSE of DWT]
  {\includegraphics[width=0.4\textwidth]{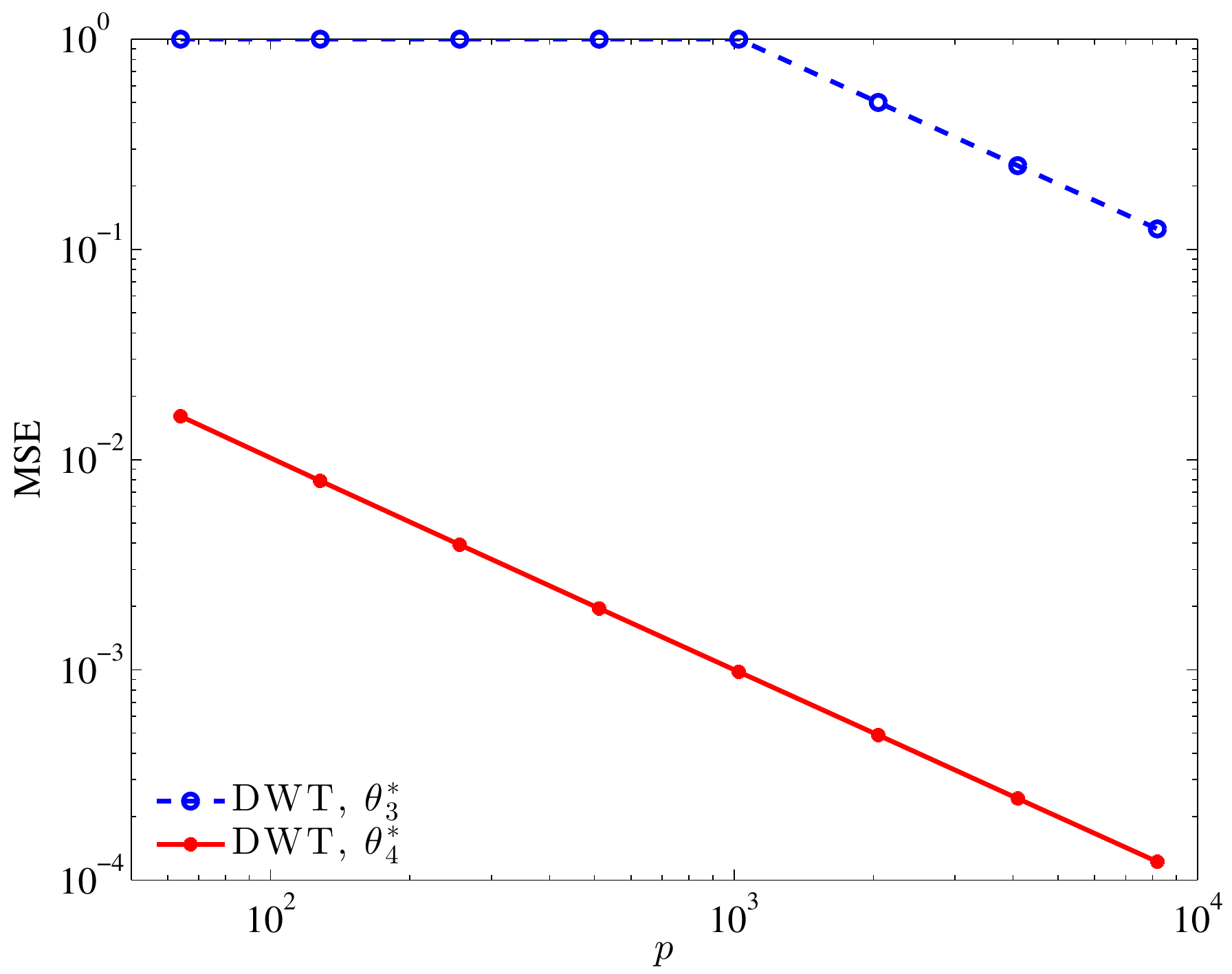} 
  \label{fig:SP_DWT}}
\caption{Signals and MSE scaling in the low-intensity regime with
  different sparse supports. $\theta^*_1$ and $\theta^*_3$ are chosen to
  make $f^*$ approximate a $\delta$-function (and thus have a large
  $\ell_2$ norm), while sparse supports 2 and 4 are randomly
  chosen. None of these signals can be rescaled to have larger
  $\ell_2$ norms without violating the nonnegativity
  constraints. $f^*_i = D\theta^*_i, i=1,\ldots,4$ are all in $\FunClass$; that is, they are
  nonnegative with unit $\ell_1$ norm. Signals like $f^*_2$ and $f^*_4$ control
  the rates in the lower bound, while signals like $f^*_1$ and $f^*_3$ control the
  rates in the upper bound. All the signals used in the experiment are $10$-sparse.}
\label{fig:sparsitypattern}
\end{figure}

Figure~\ref{fig:sparsitypattern} shows MSE behavior in the
low-intensity regime for both these signal types in the DCT and DWT bases. 
All the signals used in the simulation
belong to the function class $\FunClass$, but the signals' sparse
supports and energies vary. For
Figures~\ref{fig:DCT_theta},~\ref{fig:DCT_f} and~\ref{fig:SP_DCT} the
signals are sparse in the DCT basis, and for
Figures~\ref{fig:DWT_theta},~\ref{fig:DWT_f} and~\ref{fig:SP_DWT} the
signals are sparse in the DWT basis.  In the plots, signals
$\theta^*_1$ and $\theta^*_3$ are chosen to make $f^*_i=D\theta^*_i$
for $i = 1,3$ approximate a $\delta$-function (and thus have a large
$\ell_2$ norm). Signals $\theta^*_2$ and $\theta^*_4$ have randomly
selected supports and non-zero coefficient magnitudes equal to
$1/(p\lambda_s)$, as described above. None of these signals can be
rescaled to have larger $\ell_2$ norms without violating the
nonnegativity constraints. Signals $\theta^*_2$ and $\theta^*_4$
control the rates in the lower bound, while signals $\theta^*_1$ and
$\theta^*_3$ control the rates in the upper bound.  Details on how the
different signals were generated are in
Appendix~\ref{sec:SparseSupportDetail}.

\subsubsection{MSE vs.\ $\inty$}
\label{sec:simI}
As discussed before, the MSE's behavior changes under different
intensities. In the high-intensity setting, theoretical bounds predict
that MSE will be proportional to $1/\inty$. In the low-intensity
setting, however, $\inty$ is not a dominating factor in the lower
bound -- in fact, our bound is independent of $T$ when $T$ is below a
critical, sparsifying-basis dependent threshold.

\begin{figure}[h!]
\centering
\includegraphics[width=0.5\textwidth]{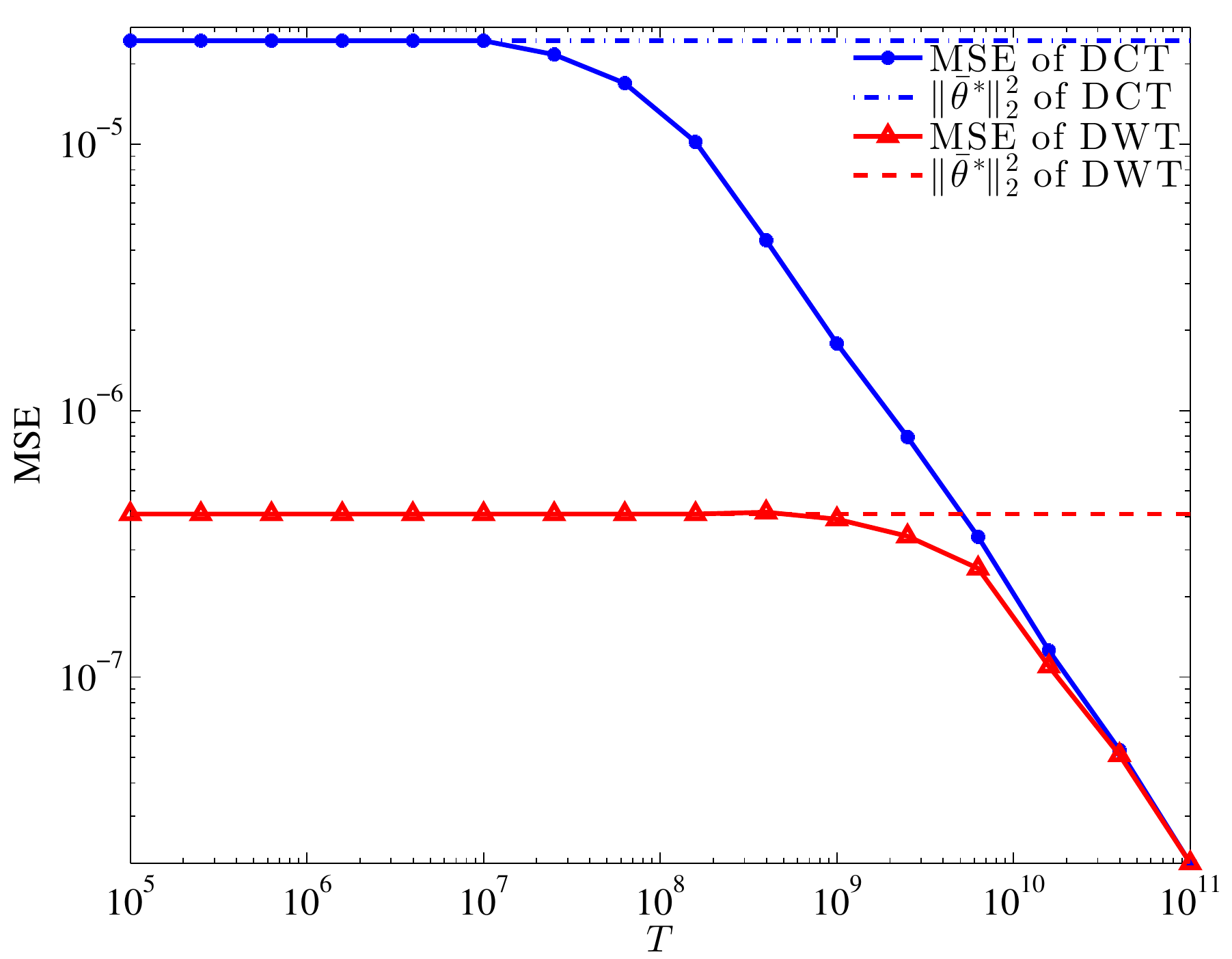}
\caption{MSE vs.\ $\inty$. The plot shows MSE of signals sparse under
  the DCT basis and of signals sparse under the DWT basis. The dashed
  lines are the values of $\|\btheta^*\|_2^2$ used in the experiments;
  this value is also an upper bound of the MSE, and we see that this
  bound is tight at low intensities. }
\label{fig:idiff}
\end{figure}

In Fig.~\ref{fig:idiff}, the plot shows how the total intensity
affects the performance. ``Elbows'' can be observed in both lines for
the DCT basis and the DWT basis, which indicate the behavior change
predicted by our theory. 
Note that in Fig.~\ref{fig:idiff} we are plotting in $\log$-$\log$
scale, and $\mbox{MSE}\propto \frac{1}{\inty} \Leftrightarrow
\log(\mbox{MSE})\propto - \log \inty$, which is exactly the linear
relationship we observe when $T$ is large. Once $\inty$ drops below
the critical value (left hand side of the plot), the MSE does not
change for different $\inty$, which reflects the bound not depending
on $\inty$ in the low-intensity regime. 
The dashed lines are the value of $\|\btheta^*\|_2^2$ used in the
experiments, which is an upper bound on the MSE. (The bound in
Thm.~\ref{thm:upper} reflects an upper bound on $\|\btheta^*\|_2^2$.)
Thus this plot suggests that MSE is proportional to $\|\btheta^*\|_2^2$ for small $T$ as suggested by our discussion.

It should also be noted that the ``elbow'' for the DCT basis arrives
at smaller $T$ than for the DWT basis, as predicted recalling that the transition between low- and high-intensity regimes occurs when $T \propto \lambda^2(D)$, and $\lambda(\bD^{{\rm    DWT}})>\lambda(\bD^{{\rm DCT}})$ when $s < \sqrt{\nc}$.

\subsubsection{MSE vs. $\nm$}
\label{sec:howmany}
One practical question that arises in many Poisson inverse problems is
how to best trade off between the number of measurements and the
amount of time spent collecting each measurement (\ie the SNR of each
measurement). Is it better to have a lot of noisy measurements, or to
have a small number of high SNR measurements?  The bounds derived in
this paper help address this question.

 In particular, Theorems~\ref{thm:lower} and~\ref{thm:upper} suggest that
the upper and lower bounds are independent of $\nm$ as long as $\nm$
is large enough to ensure that the sensing matrix $A$ satisfies
Assumptions~\ref{as:Arip1} and~\ref{as:Arip2}, an observation
conjectured in~\cite{dense_pcs}. This is contrary to most compressed
sensing results, where the error usually goes to zero as the number of
measurements increases. This interesting effect (or non-effect) of
$\nm$ is a direct result of the flux-preserving constraint, which
enforces that as $n$ increases, the number of events detected per
sensor decreases. In other words, when $\inty$ is held as a constant,
the increase of $\nm$ does not increase the overall signal-to-noise
ratio. Instead, the energy is spread out onto more observations, which
may bring more information about the signal, but also causes the
noise-per-observation to rise. 

This result suggests that  once $n$ is sufficiently large to ensure
Assumptions~\ref{as:Arip1} and~\ref{as:Arip2} are satisfied, there is no
advantage to increasing $n$ further. In other words, a small number of
high SNR measurement is sufficient, and often better in terms of
hardware costs and reconstruction computational complexity. 
This result is similar to
the finding recently reported in literature exploring the effects of
quantization on compressed sensing measurements
\cite{laska2012regime,plan2012robust}.

\begin{figure}[h!]
\centering
\includegraphics[width=0.45\textwidth]{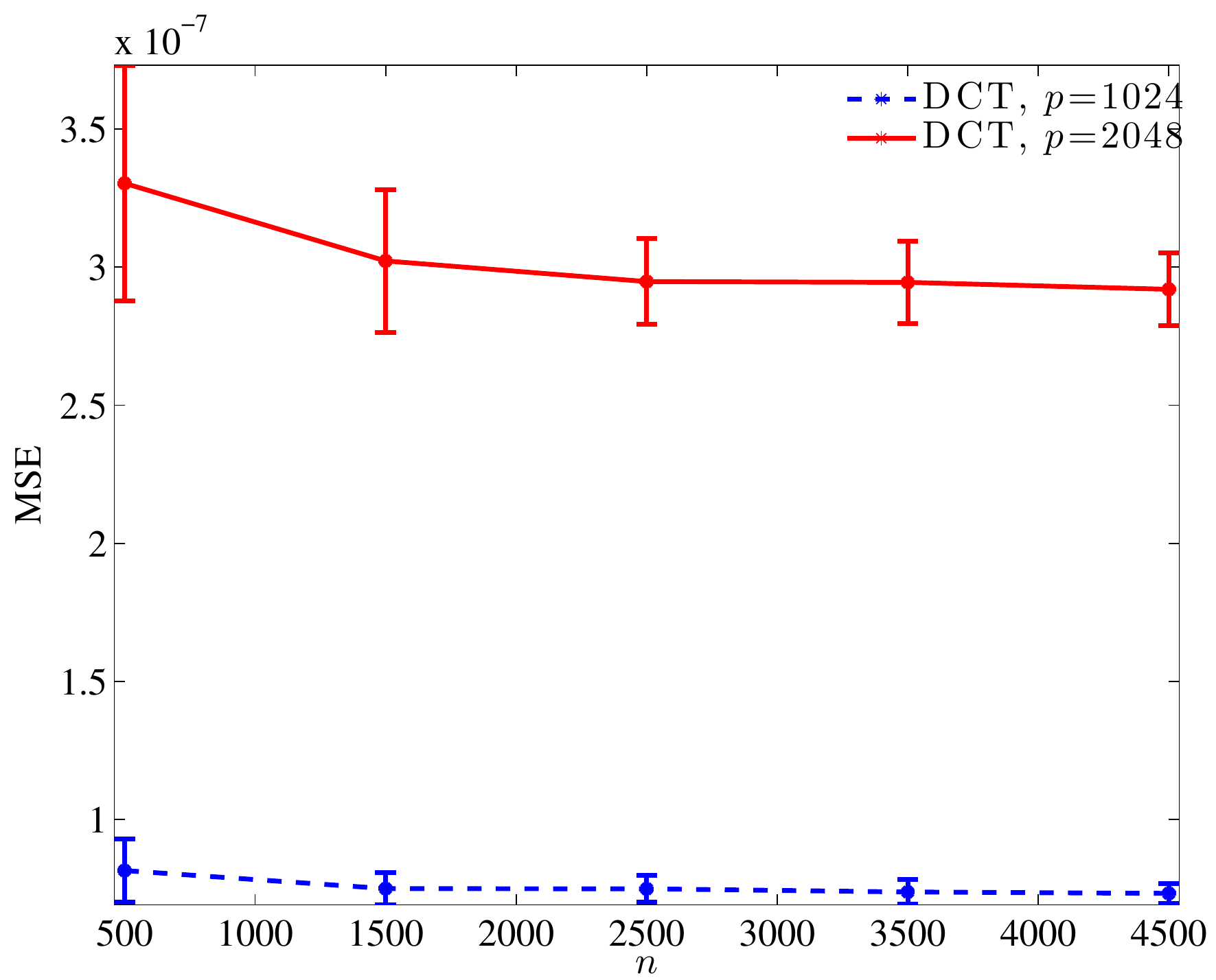}~\includegraphics[width=0.45\textwidth]{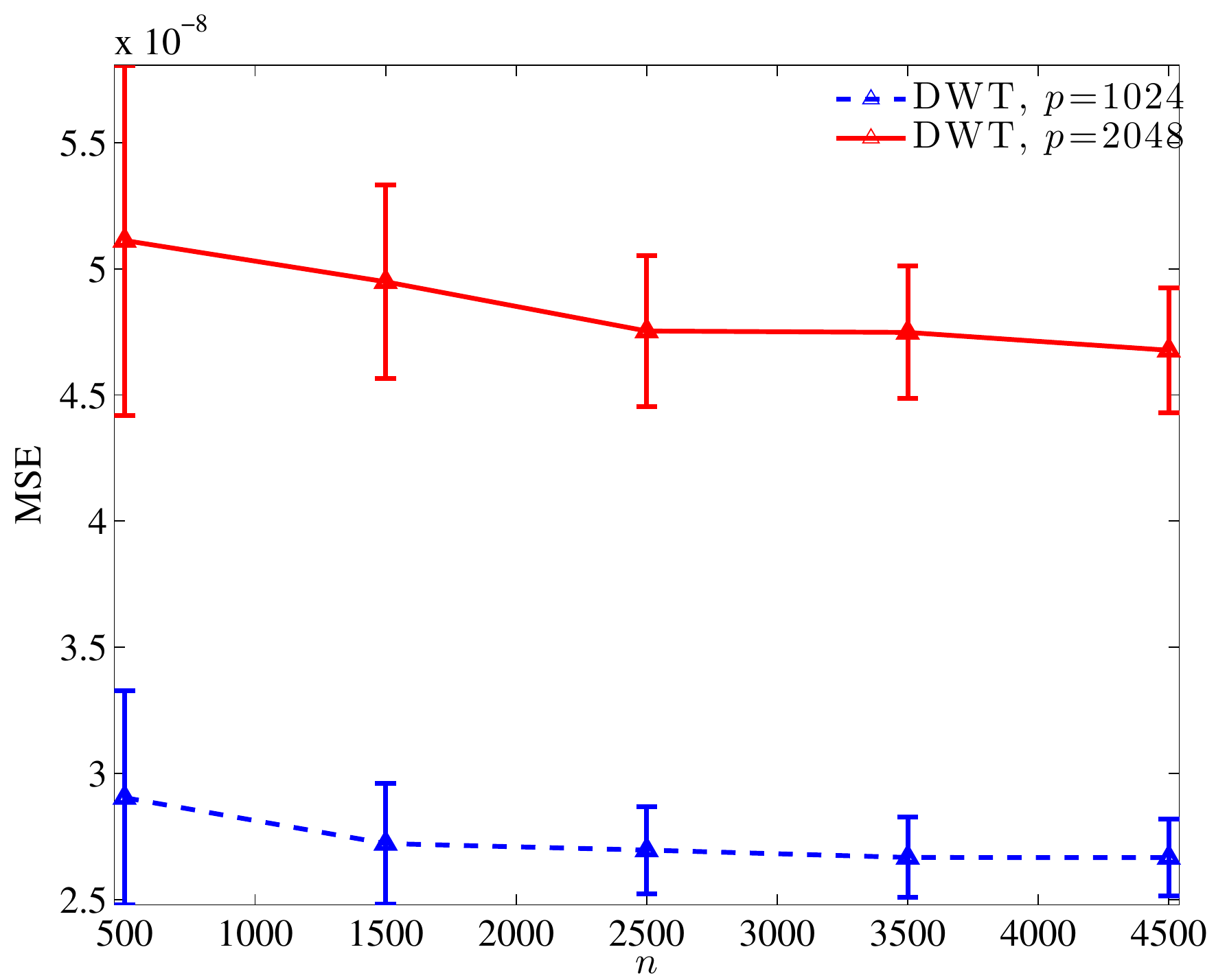}
\caption{MSE vs.\ $\nm$ ($\inty = 10^{10}, s=10$). 
The reconstruction error keeps almost constant as
$\nm$ varies within the error bars, which is consistent with the
bounds where $\nm$ has little influence.}
\label{fig:nconstant}
\end{figure}

These theoretical findings are illustrated in the below simulation.
Fig.~\ref{fig:nconstant} displays MSE of the reconstructions as a
function of $\nm$. The reconstruction error keeps almost constant as
$\nm$ varies within the error bars, which is consistent with the
bounds where $\nm$ has little influence.  

\subsubsection{MSE vs.\ $s$}

\begin{figure}[h!]
\centering
\subfloat[DCT sparsity]{\includegraphics[width=0.48\textwidth]{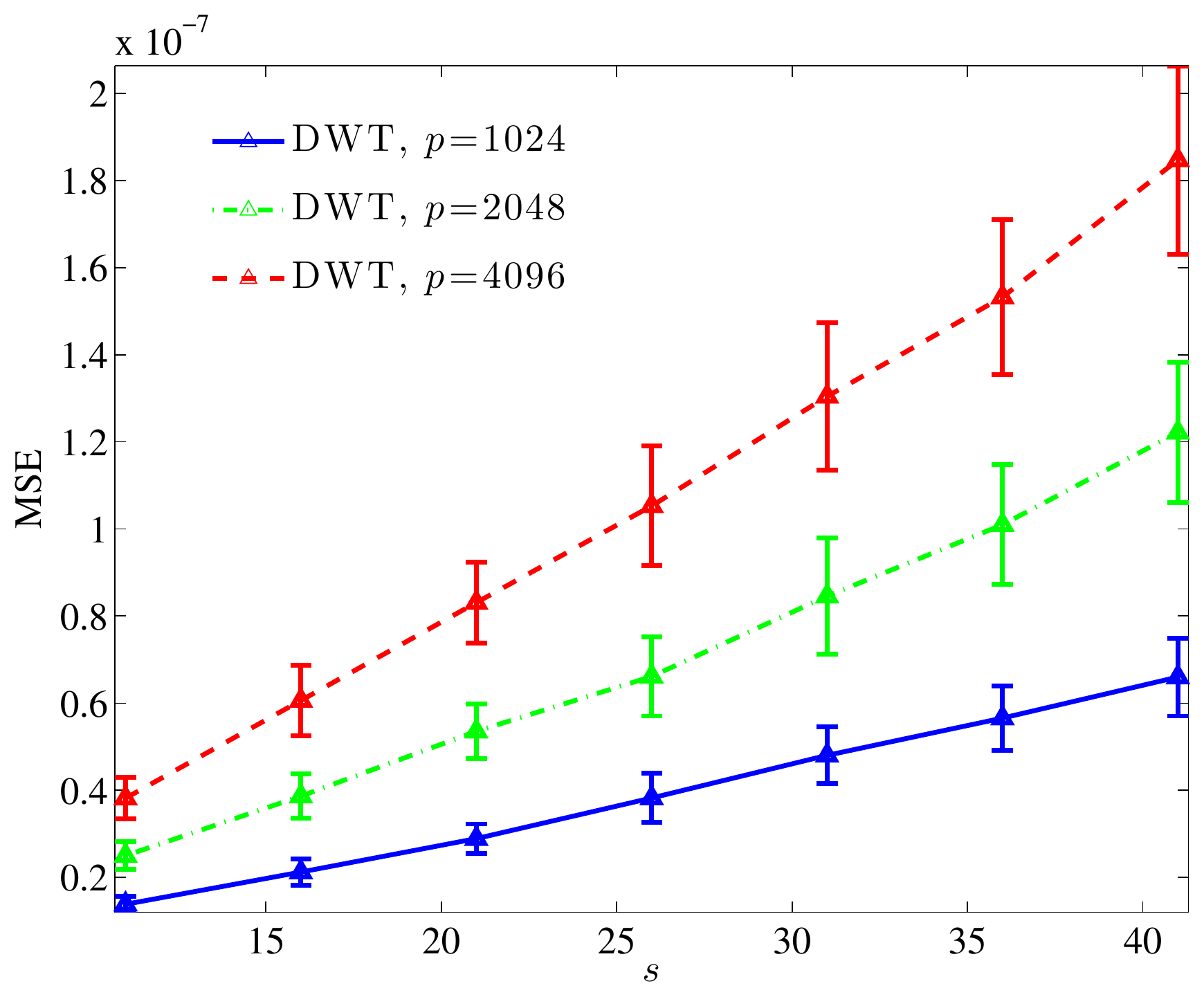}} ~ 
\subfloat[DWT sparsity]{\includegraphics[width=0.48\textwidth]{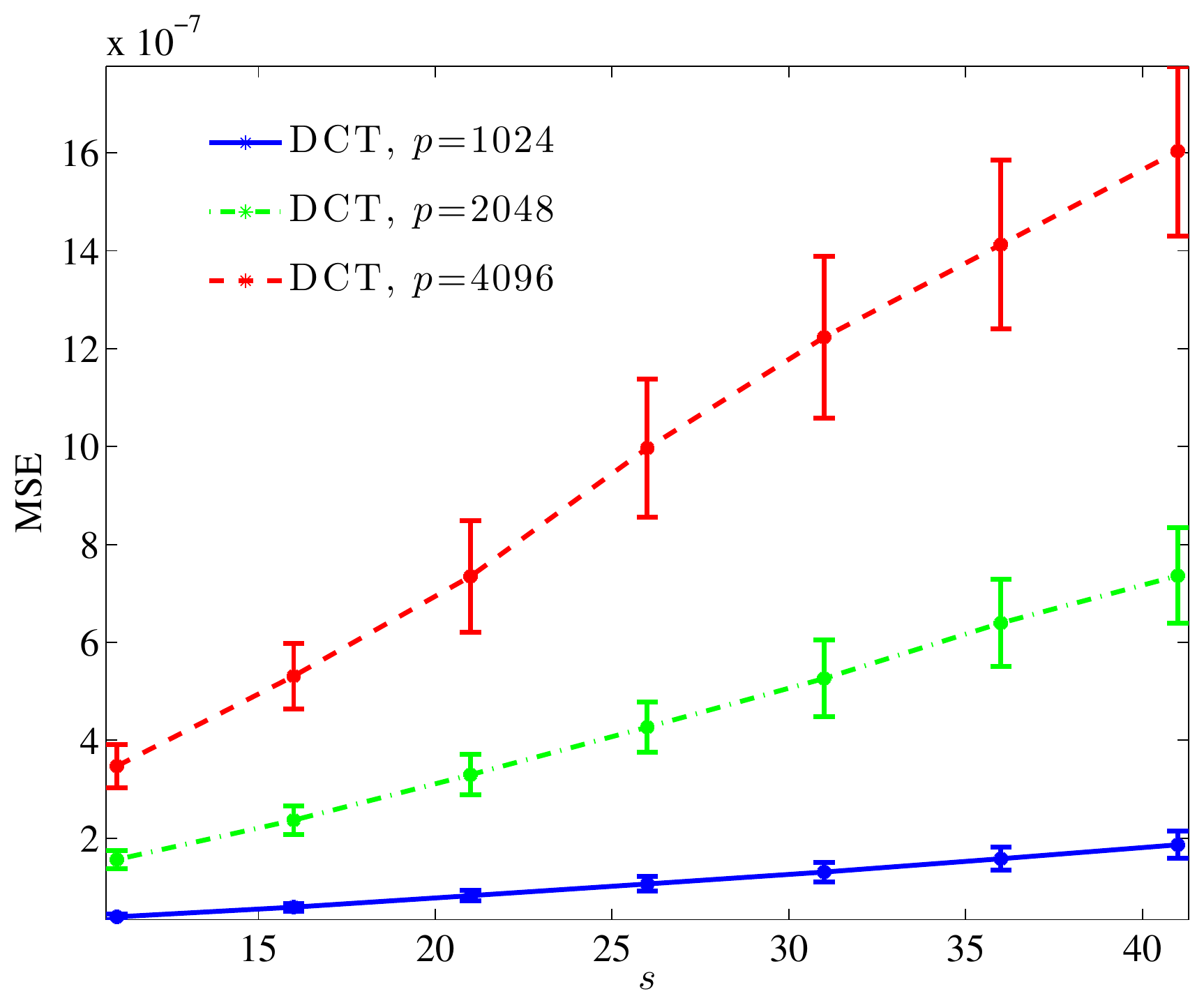}} 
\caption{MSE vs.\ $s$ in the high-intensity setting ($T = 10^{12}$ and
  $n = 500$). With both bases, the MSE grows
linearly with $s$ across a range of different $p$, which is consistent with the theoretical results.}
\label{fig:linears}
\end{figure}

Our bounds in Thms.~\ref{thm:lower} and~\ref{thm:upper} suggest the
MSE scales linearly with $s$ at high intensities (\ie when $T$ is
large). In Fig.~\ref{fig:linears}, MSE in the high-intensity ($T = 10^{12}$) is
shown for both the DCT and DWT bases for $n = 500$. With both bases, the MSE grows
linearly with $s$ across a range of different $p$, which is consistent with the theoretical results.

\subsection{Comparison of compressed sensing to downsampling}
\label{sec:dnsamp}
When designing imaging systems, one might choose between a compressed
sensing setup with incoherent projections and a ``downsampling'' setup
where we directly measure a low-resolution version of our signal (this
is specified below). The rates in Section~\ref{thm:lower}
and~\ref{thm:upper} are based on assuming a variant of a sensing
matrix which satisfies the RIP; hence standard CS schemes are capable
of achieving the reported rates. However, in practical settings one
may experience more success with a downsampling scheme. Why is this?

The key point is that many practical signals and images have
significant low-resolution content or are smooth. That is, they are not only sparse,
but the non-zero coefficients are very likely to correspond to
low-frequency or coarse-scale information. This structure in the
sparse support is not reflected by our analysis. When that structure
is present, however, simple downsampling schemes can significantly
outperform CS schemes, especially at low intensity levels. 

We define one na\"ive sensing method which uses the smoothness
of the signal as follows:
\begin{definition}[Downsampling method]
Let $\kappa \deq p/n$ be the downsampling factor (assumed to be
integer-valued). The sensing matrix for downsampling is
\begin{equation}
A^{{\rm DS}} = \mathbb{I}_{p/\kappa} \otimes \bI_{1\times \kappa},
\label{eq:dsmatrix}
\end{equation}
where {$\mathbb{I}_{p/\kappa}$} is a $p/\kappa \times p/\kappa$
identity matrix, and $\otimes$ is the Kronecker product. (That is, the
first row of $A^{\rm DS}$ is the sum of the first $\kappa$ rows of
$~\mathbb{I}_p$, the second row of $A^{\rm DS}$ is the sum of the
second $\kappa$ rows of $~\mathbb{I}_p$, etc. 
The downsampling estimator for $\theta^* \deq D^T f^*$ is 
\begin{align}
\wh{\theta}^{{\rm DS}}=\begin{bmatrix}
\frac{1}{\sqrt{p}}\\ \\
\frac{1}{\kappa\inty} \bar{D}^T A^T y^{{\rm DS}}
\end{bmatrix}.
\label{eq:dsrecover}
\end{align} 
\end{definition}

This downsampling method enhances the SNR of any individual
observation by sacrificing detailed information in the signal. Note
for the general setting of compressed sensing method, knowing that the
signal is smooth does little, if anything, to help the
reconstruction. The compressed sensing method does not distinguish
between spatial frequencies, so the reconstruction of each coefficient
is equally hard. The downsampling method makes recovering
low-frequency components significantly easier than recovering others,
and we can
obtain better performance than compressed sensing when the signal is
smooth.  

The theoretical bounds for
the downsampling method and the details of the simulation setup can
be found in Appendix~\ref{sec:dscs_detail}. In Figure~\ref{fig:dscs} we compare the performance of the
downsampling method and the compressed sensing method both
theoretically and experimentally. In the plots, $s'\le s$ reflects the
smoothness of the signal (\ie the number of low-frequency or
coarse-scale non-zero coefficients). The results show that in the
low-intensity regime, when the true signal is smooth, the
downsampling method can achieve better performance than the
compressed sensing method. As predicted by our theory, however, once
the intensity exceeds a critical threshold, we can recover both high-
and low-frequency nonzero coefficients and we begin to see an
improvement of the compressed sensing method over the downsampling method.

\begin{figure}[htpb]
\centering
\subfloat[Theoretical rates]{\includegraphics[width=0.48\textwidth]{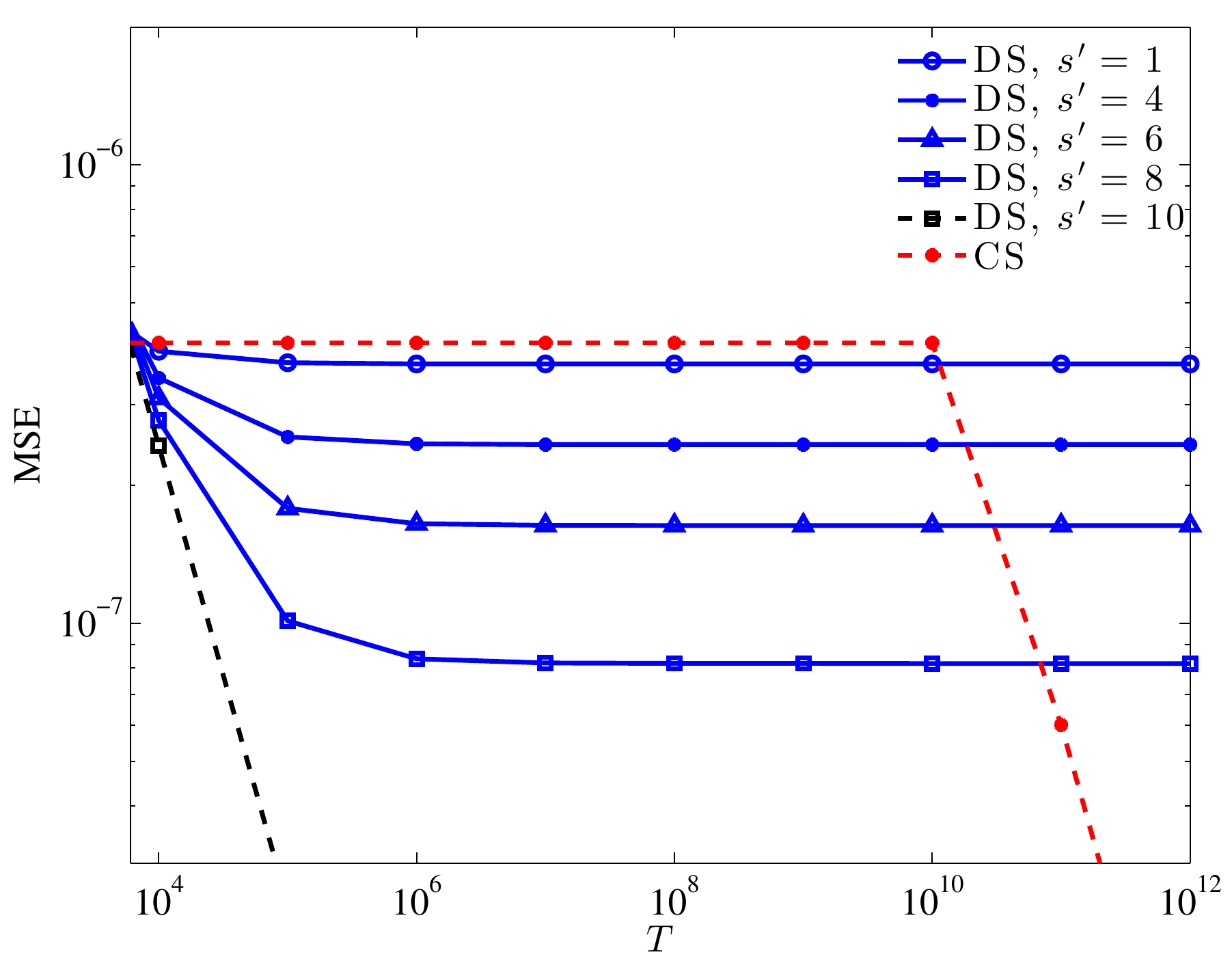}\label{fig:dscs_th}}
\subfloat[Empirical rates]{\includegraphics[width=0.48\textwidth]{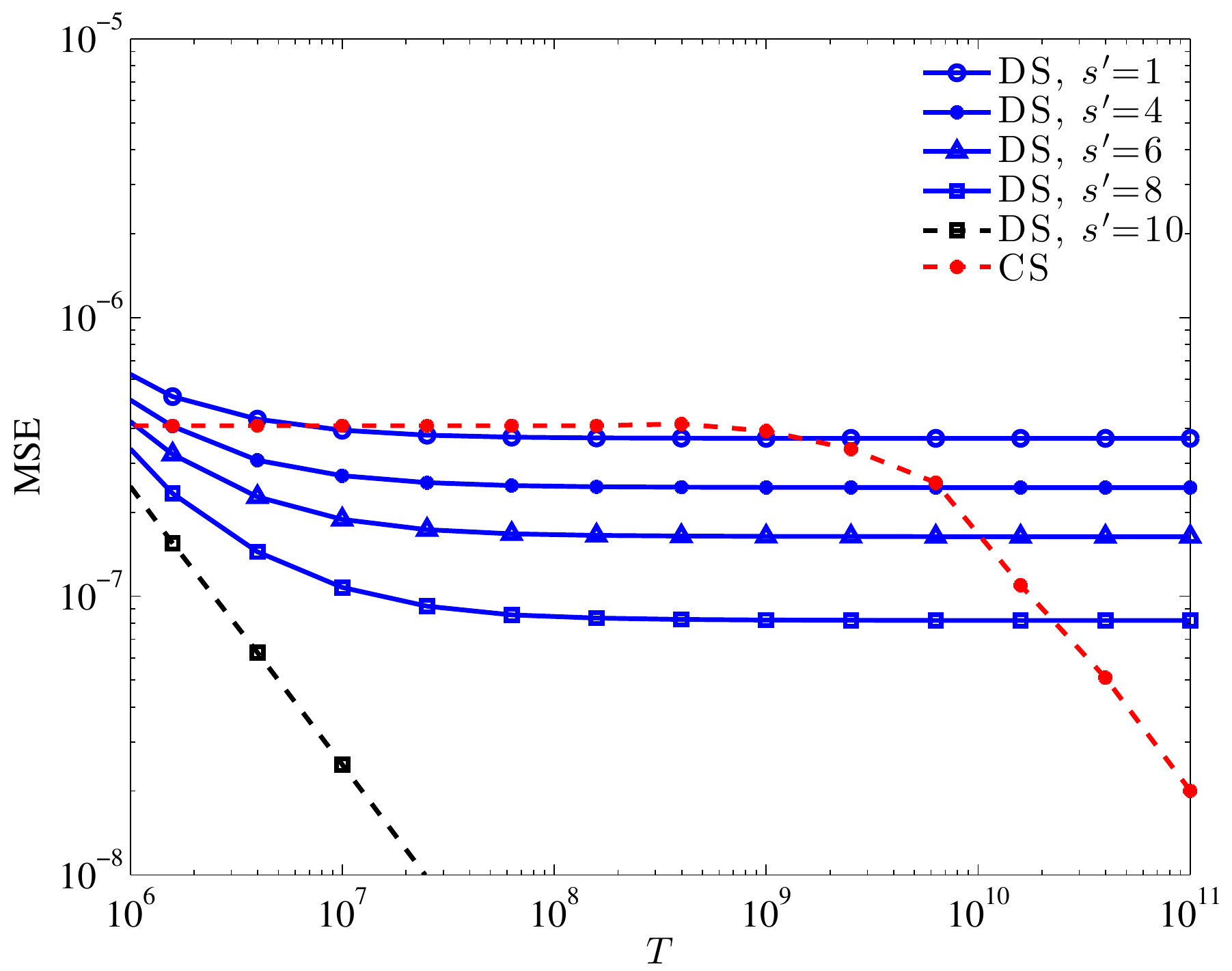}\label{fig:dscs_pr}}
\caption{Theoretical and empirical rates of downsampling and
  compressed sensing methods. $\nc=2048, K=4, \nm=\nc/K$, and $s=10$. 
  $s'$ is the
  number of coarse-scale nonzero coefficients which are directly
  measured by the proposed downsampling scheme. We see that at
  low-intensities, downsampling can yield much lower MSEs, but after
  the intensity exceeds a critical threshold, compressed sensing
  methods are able to estimate all nonzero coefficients accurately and
the MSE is better than for downsampling schemes. This effect is
predicted by our theory.}
\label{fig:dscs}
\end{figure}

\subsection{Related work}

Note that this model is different from a generalized linear model of the form
\begin{equation}
y \sim {\rm Poisson}(e^{Af^*}),
\label{eq:glm}
\end{equation}
where the exponentiation is considered to be element-wise. The model in \eqref{eq:glm} is appropriate for some problems and has been analyzed elsewhere in the CS and LASSO literature (\cf \cite{van2008high, rish2014sparse, kakade2009learning}), but is not a good representation in the motivating applications described in the introduction. For instance, in the context of imaging, the pixel intensities in $f^*$ are linearly combined by the physics of the imaging system, and that linear structure is not preserved by the exponentiation in \eqref{eq:glm}. Furthermore, in the context of \eqref{eq:model} and the above motivating applications, we face physical constraints on $A$ and $f^*$ which are not necessary in the GLM framework of \eqref{eq:glm}.

In previous performance analyses of Poisson compressed sensing, \cite{rish2009sparse} provides upper bounds of reconstruction performance for a constrained $\ell_1$-regularized estimator based on the restricted isometry property (RIP), \cite{kakade2009learning} provides upper bounds for $\ell_1$-regularized maximum likelihood estimators based on restricted eigenvalue (RE) conditions, and \cite{jia2010LASSO} focuses on understanding the impact of model mismatch and heteroscedasticity within standard LASSO algorithm, but none of these account the physical constraints on the sensing systems or provides lower bounds.
\cite{expander_pcs, dense_pcs} consider Poisson compressed sensing with similar physical constraints to this paper, but only provided upper bounds on the risk. Furthermore, in \cite{expander_pcs}, the sensing matrix is based on expander graphs, yielding different bounding techniques and algorithms. In this paper we prove both upper and lower bounds and show that they are tight in high-intensity settings. 

\section{Proofs}
\label{sec:Proof}

In this section we provide the proofs of Theorems~\ref{thm:lower} and ~\ref{thm:upper}. 

\subsection{Proof of Theorem~\ref{thm:lower}}
\label{sec:proofLower}
This proof of Theorem~\ref{thm:lower} roughly follows standard techniques for proving lower bounds on minimax rates for sparse high-dimensional problems developed in ~\cite{garvesh}. The proof involves constructing a packing set for $\FunClass$ and then applying the generalized Fano method to the packing set (see e.g.~\cite{Han94, IbrHas81, YanBar99} for details). The main challenge and novelty in the proof is adapting these standard arguments to the Poisson setting and constructing a packing set that satisfies all the physical constraints imposed in $\FunClass$ which results in the term $\frac{k}{p^2 \lambda_k ^2}$. 

\begin{proof}
We first introduce the packing sets that will be used in the proof. 
For $k=1,\ldots, s$, let
\begin{equation}
\begin{aligned}
\bar{\cH}_k \deq  \left\{\beta\in\{-1,0,+1\}^{\nc-1}:\|\beta\|_0={k}\right\}.
\end{aligned}
\notag
\end{equation}
By \cite[Lemma 4]{garvesh}, there exists a subset $\wt{\cH}_k\subseteq\bar{\cH}_k$ with cardinality $|\wt{\cH}_k|\ge\exp{(\frac{{k}}{2}\log{\frac{\nc-{k}-1}{{k}/2}})}$ such that the Hamming distance $\rho_{H}(\beta,\beta')\ge\frac{{k}}{2}$ for all $\beta,\beta'\in\wt{\cH}_k$. For a given $\alpha_k \in \reals_+$, we rescale the elements of $\wt{\cH}_k$ by $\alpha_k$ and concatenate each vector in the set with an extra entry to form the new set $\cH_{k,\alpha_k}$:
\begin{equation}
\begin{aligned}
\cH_{{k},\alpha_k} \deq \left\{\theta\in\reals^{\nc}:\theta = [1/\sqrt{\nc}, \alpha_k\beta^\T ]^\T , \beta\in\wt{\cH}\right\}.
\end{aligned}
\label{eq:defofH}
\end{equation}
Define $\eta_{\alpha_k}^2\deq \frac{k}{2}{\alpha_k}^2$, then the elements of $\cH_{k,\alpha_k}$ form a $\eta_{\alpha_k}$-packing set of $\FunClass$ in the $\ell_2$ norm. The following lemma describes several useful properties of this packing set.
\begin{lemma}
\label{lem:propertyofh}
For any $k = 1,\ldots,s$, let $\lambda_{k} = \lambda_{k}(\bD)$.
Then, the packing sets $\cH_{k,\alpha_k}$ with $0<\alpha_k\le\frac{1}{\nc\lambda_k}$ have the following properties:
\begin{enumerate}
\item{} The $\ell_2$ distance between any two points $\theta$ and $\theta'$ in $\cH_{k,\alpha_k}$ is bounded:
\begin{equation}
\begin{aligned}
\eta_{\alpha_k}^2\le\|\theta-\theta'\|_2^2\le8\eta_{\alpha_k}^2.
\end{aligned}
\notag
\end{equation}

\item{} For any $\theta\in\cH_{k,\alpha_k}$, the corresponding $f=D\theta$ satisfies:
\begin{equation}
\begin{aligned}
f_i \ge 0, \;\forall\; i \in \{1,\ldots,p\} \qquad \mbox{and} \qquad \|f\|_1 = 1.
\end{aligned}
\notag
\end{equation}

\item{} The size of the packing set 
\begin{equation}
\begin{aligned}
|\cH_{k,\alpha_k}|\ge \exp{\left(\frac{k}{2}\log{\frac{\nc-k-1}{k/2}}\right)}.
\end{aligned}
\notag
\end{equation}

\end{enumerate}

\end{lemma}
The proof of this lemma is provided in Section~\ref{sec:proofPropH}.

Throughout the proof, we also use $\Phi \deq AD$. Thus $\Phi\theta = Af$ when $f = D\theta$.
The next several steps follow the  techniques developed by \cite{Han94, IbrHas81, YanBar99}; we describe all steps for completeness. Let $M_k \deq|\cH_{k,\alpha_k}|$ be the cardinality of $\cH_{k,\alpha_k}$, and let the elements of $\cH_{k,\alpha_k}$ be denoted $\{\theta^1,\ldots,\theta^{M_k}\}$. Define a random vector $\wt{\Theta}\in\reals^\nc$ that is drawn from a uniform distribution over the packing set, and form a multi-way hypothesis testing problem where $\wt{\theta}$ is the testing result
that takes value in the packing set. Then we can bound the minimax estimation error according to~\cite{YanBar99}:
\begin{equation}
\begin{aligned}
\min_{f}\max_{f^*\in\FunClass} \expect\|f-f^*\|_2^2 \ge 
\frac{1}{4}\eta_{\alpha_k}^2\min_{\wt{\theta}}\prob[\wt{\theta}\neq\wt{\Theta}].
\end{aligned}
\label{eq:minimaxbound}
\end{equation}
By Fano's inequality and the convexity of mutual information, we have 
\begin{equation}
\begin{aligned}
\prob[\wt{\theta}\neq\wt{\Theta}]
&\ge 1-\frac{I(y;\wt{\Theta}) + \log{2}}{\log{M_k}},
\end{aligned}
\label{eq:pbound}
\end{equation}
where $y \sim \mbox{Poisson}(\inty\Phi\wt{\Theta})$ is the Poisson observation, 
$I(y;\wt{\Theta})$ is the mutual information between random variables $y$ and $\wt{\Theta}$. Additionally, the mutual information can be bounded by the average of the K-L divergence between $p(y|\inty\Phi\si)$ and $p(y|\inty\Phi\sj)$ for all $\si, \sj \in \cH_{k,\alpha_k}$,\ie
\begin{equation}
\begin{aligned}
I(y;\wt{\Theta})
&\le \frac{1}{\binom{M_k}{2}}\sum_{\substack{i,j=1,\ldots,M_k \\ i\neq j}}\KL(p(y|\inty\Phi\si)\|p(y|\inty\Phi\sj)),
\end{aligned}
\label{eq:ibound}
\end{equation}
according to ~\cite{Han94}.

The following lemma provides an upper bound for the K-L divergence of Poisson distributions in terms of the squared $\ell_2$-distance.

\begin{lemma}
\label{lem:PoisKL}
Let $p(y|\mu)$ denote the multivariate Poisson distribution with mean parameter $\mu\in\reals^{\nm}_+$. For $\mu_1,\mu_2\in\reals_+^{\nm}$, if there exists some value $c>0$ such that $\mu_2\succeq c\bI_{\nm\times 1}$, then the following holds:
\begin{equation}
\KL(p(y|\mu_1)\|p(y|\mu_2)) \le \frac{1}{c}\|\mu_1-\mu_2\|_2^2.
\notag
\end{equation}
\end{lemma}

The proof can be found in Section~\ref{sec:proofofKL}. Assumption~\ref{as:wtA} ensures $Af^*$ is bounded, as we have the following lemma:
\begin{lemma}
\label{lem:AfBound}
If the sensing matrix $A$ satisfies Assumption~\ref{as:wtA}, then for all nonnegative $f$ s.t.  $\|f\|_1=1$, we have
\begin{equation}
\eps\bI_{\nm\times 1} \preceq Af \preceq \frac{1}{\nm}\bI_{\nm\times 1}.
\label{eq:prop2}
\end{equation} 
\end{lemma}
The proof is provided in Section~\ref{sec:proofofAbound}.

By Lemmas~\ref{lem:propertyofh} and~\ref{lem:AfBound}, we have $\Phi\sj \succeq {\eps}\bI_{\nc\times 1}$, \ie $\inty\Phi\sj \succeq {\frac{\inty}{2\nm}}\bI_{\nc\times 1}$. Thus by Lemma~\ref{lem:PoisKL},
\begin{equation}
\begin{aligned}
KL(p(y|\inty\Phi\si)\|p(y|\inty\Phi\sj))
&\le \frac{2\nm}{\inty}\|\inty\Phi(\si-\sj)\|_2^2 \\
&= {2\nm\inty}\|\Phi(\si-\sj)\|_2^2.
\end{aligned}
\notag
\end{equation}
Let $f^i \deq D\si$, and $f^j \deq D\sj$, then
\begin{equation}
\begin{aligned}
\|\Phi(\si-\sj)\|_2^2 
&= \|A(f^i-f^j)\|_2^2\\
&= \left\|\frac{1}{\resc}\wt{A}(f^i-f^j) - 
	\frac{\shft}{\resc}\bI_{\nm\times\nc}(f^i-f^j)\right\|_2^2 \\
&= \left\|\frac{1}{\resc}\wt{A}D(\si-\sj) \right\|_2^2 \\
&\le \frac{\Ru}{4(a_u-a_\ell)^2\nm} \|(\si-\sj)\|_2^2.
\end{aligned}
\notag
\end{equation}
where the last equality is a result of $\|f^i\|_1 = \|f^j\|_1$ which leads to $\bI_{\nm\times\nc}(f^i-f^j)=0$, and the inequality follows from Assumption~\ref{as:Arip1}.
The construction of our packing set ensures that $\|\si-\sj\|_2^2\le 8\eta_{\alpha_k}^2$, so
\begin{equation}
\begin{aligned}
\|\Phi(\si-\sj)\|_2^2 
&\le \frac{2\Ru}{(a_u-a_\ell)^2\nm}\eta_{\alpha_k}^2.
\end{aligned}
\notag
\end{equation}
As a result,
\begin{equation}
\begin{aligned}
KL(p(y|\inty\Phi\si)\|p(y|\inty\Phi\sj))
&\le {2\nm\inty}\|\Phi(\si-\sj)\|_2^2 \\
&\le \frac{4\Ru\inty}{(a_u-a_\ell)^2}\eta_{\alpha_k}^2.
\end{aligned}
\label{eq:klbound}
\end{equation}
By Ineq.~[\eqref{eq:ibound}, \eqref{eq:klbound}], we can bound the mutual information
\begin{equation}
\begin{aligned}
I(y;\wt{\Theta})
&\le \frac{1}{\binom{M}{2}}\sum_{i\neq j}KL(p(y|\inty\Phi\si)\|p(y|\inty\Phi\sj)) \\
&\le \max_{i\neq j}KL(p(y|\inty\Phi\si)\|p(y|\inty\Phi\sj)) \\
&\le \frac{4\Ru\inty}{(a_u-a_\ell)^2}\eta_{\alpha_k}^2.
\end{aligned}
\label{eq:ibound_f}
\end{equation}
Then, by Ineq.~[\eqref{eq:pbound}, \eqref{eq:ibound_f}] we can bound the probability of classification error by
\begin{equation}
\begin{aligned}
\prob[\wt{\theta}\neq\wt{\Theta}]
&\ge 1-\frac{I(y;\wt{\Theta}) + \log{2}}{\log{M_k}} \\
&\ge 1-\frac{\frac{4\Ru\inty}{(a_u-a_\ell)^2}\eta_{\alpha_k}^2 + \log{2}}{\frac{k}{2}\log{\frac{\nc-k-1}{k/2}}}.
\end{aligned}
\label{pbound_f}
\end{equation}
We need to ensure this probability is bounded below by a positive constant; in the below we use the constant $1/4$. This constant is guaranteed as long as we have 
\begin{equation}
{\frac{k}{2}\log{\frac{\nc-k-1}{k/2}}} \ge 2\log2
\label{eq:Mbound1}
\end{equation} 
and 
\begin{equation}
{\frac{k}{2}\log{\frac{\nc-k-1}{k/2}}}\ge\frac{16\Ru\inty}{(a_u-a_\ell)^2}\eta_{\alpha_k}^2,
\label{eq:Mbound2}
\end{equation} 
because then we can bound
\begin{equation}
\begin{aligned}
\prob[\wt{\theta}\neq\wt{\Theta}]
&\ge1 - \frac{\frac{4\Ru\inty}{(a_u-a_\ell)^2}\eta_{\alpha_k}^2}
{\frac{k}{2}\log{\frac{\nc-k-1}{k/2}}} 
- \frac{\log{2}}{\frac{k}{2}\log{\frac{\nc-k-1}{k/2}}} \\
& \ge 1 - \frac{1}{4} - \frac{1}{2} \\
& \ge \frac{1}{4}.
\end{aligned}
\notag
\end{equation}
First we verify Ineq.~\eqref{eq:Mbound1}. For $k\ge2$, 
\begin{equation}
\begin{aligned}
\frac{k}{2}\log{\frac{\nc-k-1}{k/2}} &\ge \log{\frac{\nc-k-1}{k/2}} \\
&\ge \log{\frac{\nc-s-1}{s/2}}\\
&> \log{\frac{\nc-(s+1)}{(s+1)/2}} \\
&\ge \log{\frac{2(s+1)}{s/2}} \\
& = 2\log2,
\end{aligned}
\notag
\end{equation}
where the last inequality is a result of $\nc\ge3(s+1)$. For $k=1$, 
\begin{equation}
\begin{aligned}
\frac{k}{2}\log{\frac{\nc-k-1}{k/2}} &= \frac{1}{2}\log{\frac{\nc-2}{1/2}} \\
&\ge \frac{1}{2}\log{\frac{10-2}{1/2}} \\
& =  2\log2,
\end{aligned}
\notag
\end{equation}
where the inequality is a result of $\nc\ge 10$. Now to ensure Ineq.~\eqref{eq:Mbound2}, we need
\begin{equation}
\begin{aligned}
\frac{\Ru\inty}{(a_u-a_\ell)^2}\eta_{\alpha_k}^2
\le\frac{k}{32}\log{\frac{\nc-k-1}{k/2}},
\end{aligned}
\notag
\end{equation}
which leads to
\begin{equation}
\begin{aligned}
\eta_{\alpha_k}^2 &\le \frac{(a_u-a_\ell)^2k}{32\Ru\inty}\log{\frac{\nc-k-1}{k/2}}.
\end{aligned}
\notag
\end{equation}
Next recall the condition of Lemma~\ref{lem:propertyofh} that $0<\alpha_k\le\frac{1}{\nc\lambda_k}$. This combined with the above suggests $\alpha_k$ must be chosen so that
\begin{equation}
\begin{aligned}
\eta_{\alpha_k}^2 = \min\left\{ 
\frac{k}{2}\left(\frac{1}{\nc\lambda_k}\right)^2, 
\frac{(a_u-a_\ell)^2k}{32\Ru\inty}\log{\frac{\nc-k-1}{k/2}}
\right\},
\end{aligned}
\notag
\end{equation}
which is satisfied by our choice of $\alpha_k$. A lower bound is then proved using the packing set $\cH_{\alpha_k}$ such that
\begin{equation}
\begin{aligned}
\min_{f}\max_{f^*\in\FunClass} \expect\|f-f^*\|_2^2
&\ge\frac{1}{4}\eta_{\alpha_k}^2\prob[\wt{\theta}\neq\wt{\Theta}] \geq \frac{1}{16} \eta_{\alpha_k}^2\\
&\ge \min \left\{ \frac{1}{32} \frac{k}{\nc^2\lambda_k^2}, \frac{(a_u-a_\ell)^2k}{512\Ru\inty}\log{\frac{\nc-k-1}{k/2}}\right\}.
\end{aligned}
\notag
\end{equation}
For each level $k = 1,\ldots,s$, the above is a separate lower bound for the minimax risk. To make the bound as tight as possible, we take the maximum of the bounds for all $s$ sparsity levels, which results in
\begin{equation}
\begin{aligned}
\min_{f}\max_{f^*\in\FunClass} \expect\|f-f^*\|_2^2
\ge \max_{1\le k\le s} \left\{ \min \left( \frac{1}{32} \frac{k}{\nc^2\lambda_k^2}, \frac{(a_u-a_\ell)^2k}{512\Ru\inty}\log{\frac{\nc-k-1}{k/2}} \right) \right\}.
\end{aligned}
\notag
\end{equation}
Then there exists some constant $C_L>0$ such that
\begin{equation}
\begin{aligned}
\min_{f}\max_{f^*\in\FunClass} \expect\|f-f^*\|_2^2
\ge C_L\max_{1\le k\le s} \left\{ \min \left( \frac{k}{\nc^2\lambda_k^2}, \frac{k}{\Ru\inty}\log{\frac{\nc-k-1}{k/2}} \right) \right\},
\end{aligned}
\notag
\end{equation}
which completes the proof.
\end{proof}

\subsubsection{Proof of Lemma~\ref{lem:propertyofh}}
\label{sec:proofPropH}
\begin{proof} We first prove property 1. Let $\btheta = [\theta_2,\ldots,\theta_{\nc}]^\T, \btheta' = [\theta'_2,\ldots,\theta'_{\nc}]^\T$. Note that the first elements of $\theta$ and $\theta'$ are equal, \ie $\theta_1 = \theta'_1$, thus $\|\theta - \theta'\|_2^2 = \|\btheta - \btheta'\|_2^2$. By construction, the Hamming distance between $\theta$ and $\theta'$ is at least ${k}/2$. Because $\btheta,\btheta'\in\{0,\pm \alpha_k\}^{\nc-1}$, we have
\begin{equation}
\|\btheta - \btheta'\|_2^2\ge \frac{k}{2}\alpha_k^2 = \eta_{\alpha_k}^2.
\notag
\end{equation}
Next, note that $\|\btheta\|_0 = \|\btheta'\|_0 = k$. Thus the maximum distance is
\begin{equation}
\|\btheta - \btheta'\|_2^2 \le k(\alpha_k - (-\alpha_k))^2 = 4k\alpha_k^2 = 8\eta_{\alpha_k}^2,
\notag
\end{equation}
and we have property 1 proved.

For property 2, we have that for any $\theta\in\cH_{k,\alpha_k}$, $\theta_1 = 1/\sqrt{\nc}$, and $|\theta_i| \le \alpha_k \le \frac{1}{\nc\lambda_k}, i\neq 1$. Let $\btheta = [\theta_2,\ldots,\theta_{\nc}]^\T$ and $\beta = \btheta/\alpha_k \in\{0,\pm1\}^{\nc-1}$ be the corresponding element in $\bar{\cH_{k}}$, then
\begin{equation}
\begin{aligned}
f_i	& = (D\theta)_i \\
	& = D_{i,1}\theta_1 + \sum_{j=2}^{\nc}D_{i,j}\theta_j \\
	& = \frac{1}{\sqrt{\nc}}\theta_1 + \alpha_k\langle \bD_i, \beta \rangle \\
	& \ge \frac{1}{\sqrt{\nc}}\frac{1}{\sqrt{\nc}} - \alpha_k |\langle \bD_i, \beta\rangle| \\
	& \ge \frac{1}{\nc} - \alpha_k \lambda_k.
\end{aligned}
\notag
\end{equation}
By $\alpha_k\le\frac{1}{\nc\lambda_k}$,
\begin{equation}
\begin{aligned}
f_i	& \ge \frac{1}{\nc} - \alpha_k \lambda_k \\
	& \ge \frac{1}{\nc} - \frac{1}{\nc\lambda_k}\lambda_k  \\
	& = 0.
\end{aligned}
\notag
\end{equation}
Then, by the properties of $D$: 
\begin{equation}
\begin{aligned}
\|D\theta\|_1 &= \sqrt{\nc}\theta_1 = 1.
\end{aligned}
\notag
\end{equation}

Property 3 is satisfied by construction. The proof of the size of the packing set can be found in \cite[Lemma 4]{garvesh}.
\end{proof}

\subsubsection{Proof of Lemma~\ref{lem:PoisKL}}
\label{sec:proofofKL}
\begin{proof}
The proof is straight-forward. For Poisson distributions, we have
\begin{equation}
\begin{aligned}
KL(p(y|\mu_1)\|p(y|\mu_2))
&=\sum_{i=1}^{\nm}\left[(\mu_1)_i\log{\frac{(\mu_1)_i}{(\mu_2)_i}}-(\mu_1)_i+(\mu_2)_i\right] \\
&\le \sum_{i=1}^\nm\left[(\mu_1)_i\left(\frac{(\mu_1)_i}{(\mu_2)_i}-1\right) - (\mu_1)_i+ (\mu_2)_i\right] \\
&= \sum_{i=1}^{\nm}\frac{1}{(\mu_2)_i}\left[(\mu_1)_i^2-2(\mu_1)_i(\mu_2)_i+(\mu_2)_i^2\right] \\
&\le \frac{1}{\min_{i}(\mu_2)_i}\sum_{i=1}^{\nm}\left[(\mu_1)_i^2-2(\mu_1)_i(\mu_2)_i+(\mu_2)_i^2\right] \\
&\le \frac{1}{c}\|\mu_1-\mu_2\|_2^2,
\end{aligned}
\notag
\end{equation}
where the last inequality is a result of $\mu_2\succeq c\bI_{\nm\times1}$.
\end{proof}

\subsubsection{Proof of Lemma~\ref{lem:AfBound}}
\label{sec:proofofAbound}
\begin{proof}
For any $A$ satisfying Assumption~\ref{as:wtA}, by Eq.~\ref{eq:rangeofA}, we have
$A_{i, j} \in [\frac{1}{2\nm},\frac{1}{\nm}]$.
Thus for any nonnegative $f$ satisfying  $\|f\|_1=1$, by H\"older's inequality for $i=1,\ldots,\nc$, we have
\begin{equation}
\begin{aligned}
(Af)_i = \sum_{j=1}^{\nc} A_{i,j}f_j 
\ge \|f\|_1\min_{j =1,\ldots,\nc} A_{i,j} 
\ge \eps,
\end{aligned}
\notag
\end{equation}
and
\begin{equation}
\begin{aligned}
(Af)_i = \sum_{j=1}^{\nc} A_{i,j}f_j 
\le \|f\|_1\max_{j =1,\ldots,\nc} A_{i,j} 
\le \frac{1}{\nm}.
\end{aligned}
\notag
\end{equation}
\end{proof}

\subsection{Proof of Theorem~\ref{thm:upper}}
\label{sec:proofUpper}
\begin{proof}
Here we use the proof technique developed in \cite{dense_pcs}, but
modify it to provide sharper upper bounds. First, we do not assume
$f^*$ bigger than some positive threshold
$c/\nc\bI_{\nc\times1}$. Instead, we work with a more general set of
signals  where we only assume $f^*$ is nonnegative. We also modify the
feasible set of estimators and use slightly different assumptions on
the sensing matrix $A$ to further tighten the bound. The main steps of
the proof include: 
\begin{enumerate}
\item Construct a specific finite set $\Gamma$ of quantized feasible estimators, and design a penalty function $\pen(\cdot)$;
\item Use the assumptions on $A$ and the Li-Barron theorem on the
  error of penalized likelihood estimators (see \cite{dense_pcs, li1999mixture, kolaczyk2004multiscale}) to bound
\begin{equation}
\|\wh{f} - f^*\|_2^2 \le c\min_{(D^\T f)\in\Gamma}\left\{\|{f} - f^*\|_2^2 + \frac{1}{\inty} \pen(f)\right\},
\notag
\end{equation}
for some constant $c$;
\item Choose one estimator in the feasible set, and compute its cost as the upper bound.

\end{enumerate}
Let
$L \deq \max_{j,k} |D_{j,k}|$ 
and note
\begin{equation}
\|\theta^*\|_{\infty} = \|D^\T f^*\|_\infty \le
\max_{j,k}|D_{j,k}|\|f^*\|_1 = \max_{i,j}|D_{i,j}| = L.
\notag
\end{equation}
Next define the quantization level 
\begin{equation}
\begin{aligned}
K^* &\deq \sqrt{\frac{\Ru L^2\inty\log2}{2(a_u-a_{\ell})^2}}, \\
K &\deq \left\lceil K^* \right\rceil.
\end{aligned}
\label{eq:K}
\end{equation}
Let
\begin{equation}
\begin{aligned}
\mathcal{B} \deq 
\begin{cases}
\left\{0\right\}, & K=1 \\
\left\{-L, -\frac{K-3}{K-1}L, -\frac{K-5}{K-1}L, \ldots, \frac{K-3}{K-1}L, L \right\}, &K>1,
\end{cases}
\end{aligned}
\notag
\end{equation}
where $|\mathcal{B}| = K$. We now consider the case where $K>1$ (which also implies $K^*>1$); the
case where $K=1$ will be addressed later. 
Define the following sets:
\begin{subequations}
\begin{align}
\Theta_K \deq &\left\{
 \theta\in\reals^\nc: \theta_1 = \nc^{-1/2}, \|\theta\|_0\le s+1, \theta_i \in \mathcal{B}  \mbox{ for } i = 2, \ldots, \nc \right\}, \label{eq:Theta} \\
\mathcal{F}\deq &\{f=D\theta: \theta\in\Theta_K\},\\
\mathcal{C} \deq &\left\{ g\in\reals_+^\nc:  \|g\|_1 = 1 \right\},\\
\Gamma \deq &\left\{ \wt{\theta} = D\wt{f}: f \in \mathcal{F}, \wt{f} = \argmin_{g \in \mathcal{C}} \|f-g\|_2\right\}.
\label{eq:defofGamma}
\end{align}
\end{subequations}
Further define the penalty function 
\begin{equation}
\pen(\theta) = 
2\log_2(s)+2\|{\btheta}\|_0\log_2 (\nc K).
\label{eq:defPen}
\end{equation}
This penalty function is calculated from the code length of a prefix code for all $\theta \in\Theta_K$. Since we know $\theta_1$ exactly, we only need to code $\bar{\theta} = [\theta_2,\ldots,\theta_{\nc}]^\T$.
The prefix code can be generated using the following steps:
\begin{enumerate}
\item First, we encode $\|\bar{\theta}\|_0$, the number of non-zeros of $\bar{\theta}$, which requires $\log_2 s$ bits. 
\item Second, we encode each of the locations of the $\|\bar{\theta}\|_0$ non-zero elements. Each location takes $\log_2\nc$ bits, so the overall length is $\|\bar{\theta}\|_0\log_2\nc$ bits.
\item Last, we encode the value of each non-zero element. Since each value is quantized to one of the $K$ bins, this step will take $\|\bar{\theta}\|_0\log_2K$ bits.
\end{enumerate}
Thus the prefix code length is $\log_2(s)+\|{\btheta}\|_0\log_2 (\nc K)$. Note our penalty function is twice the code length. This scaling is needed later for the proof of Eq.~\eqref{eq:Becca26}. 
By having a corresponding prefix code, the penalty function (rescaled by $1/2$) satisfies the Kraft inequality:
$
\sum_{\theta \in\Theta_K}e^{-\pen(\theta)/2}\le 1.
$
This code can also be applied to the set $\Gamma$ by coding $\wt{\theta}$ using the same code for the corresponding $\theta\in\Theta_K$, where we also have
\begin{equation}
\begin{aligned}
\sum_{\wt{\theta} \in\Gamma}e^{-\pen(\wt{\theta})/2}\le 1.
\end{aligned}
\label{eq:Kraft}
\end{equation}

By Assumption~\ref{as:Arip2} and $\|f^*\|_1 = \|\wh{f}\|_1$, we have 
\begin{equation}
\begin{aligned}
\|f^*-\wh{f}\|_2^2 
	& \le \frac{1}{\Rl }\|\wt{A}(f^*-\wh{f})\|_2^2 \\
	& \le  \frac{(\resc)^2}{\Rl }\|A(f^*-\wh{f})\|_2^2.
\end{aligned}
\label{eq:uppera}
\end{equation}
Then,
\begin{equation}
\begin{aligned}
\|A(f^*-\wh{f})\|_2^2 
	&= \sum_{i=1}^{\nm}\left[(Af^*)_i-(A\wh{f})_i\right]^2 \\
	&= \sum_{i=1}^{\nm}\left[(Af^*)_i^{1/2}-(A\wh{f})_i^{1/2}\right]^2\left[(Af^*)_i^{1/2}+(A\wh{f})_i^{1/2}\right]^2 \\
	&\le 4\max_{i}\{(Af^*)_i, (A\wh{f})_i\}\sum_{i=1}^{\nm}\left[(Af^*)_i^{1/2}-(A\wh{f})_i^{1/2}\right]^2 \\
	&\le \frac{4}{\nm} \sum_{i=1}^{\nm}\left[(Af^*)_i^{1/2}-(A\wh{f})_i^{1/2}\right]^2,
\end{aligned}
\label{eq:upperb}
\end{equation}
where the second equality holds because of the non-negativity
constraints. The last inequality holds because $D^\T\wh{f}\in\Gamma
\Rightarrow \wh{f}\in\mathcal{C} \Rightarrow \wh{f}$ is nonnegative
and $\|\wh{f}\|_1=1$, so we can apply Lemma~\ref{lem:AfBound} to both
$\wh{f}$ and $f^*$.

Following the technique used in \cite[Eq. (25) and Ineq. (26)]{dense_pcs} (where we need the Kraft inequality in~\eqref{eq:Kraft}), we further have
\begin{equation}
\begin{aligned}
&\expect\left[\sum_{i=1}^{\nm}\left((\inty Af^*)_i^{1/2}-(\inty A\wh{f})_i^{1/2}\right)^2\right] \le \min_{(D^\T f)\in\Gamma }\left[\KL\left(p(\cdot|\inty Af^*)\|p(\cdot|\inty Af)\right) + \pen(f)\right],
\end{aligned}
\label{eq:Becca26}
\end{equation}
which is equivalent to
\begin{equation}
\begin{aligned}
&\expect\left[\sum_{i=1}^{\nm}\left((Af^*)_i^{1/2}-(A\wh{f})_i^{1/2}\right)^2\right] \le \frac{1}{\inty}\min_{(D^\T f)\in\Gamma }\left[\KL\left(p(\cdot|\inty Af^*)\|p(\cdot|\inty Af)\right) + \pen(f)\right].
\end{aligned}
\notag
\end{equation}
By Lemma~\ref{lem:PoisKL} and $Af\succeq \eps\bI_{\nm\times1}$ from Lemma~\ref{lem:AfBound}, we have
\begin{equation}
\begin{aligned}
\KL\left(p(\cdot|\inty Af^*)\|p(\cdot|\inty Af)\right) 
&\le \frac{2\nm}{\inty}\|\inty A(f^*-f)\|_2^2 = {2\nm\inty}\|A(f^*-f)\|_2^2.
\end{aligned}
\notag
\end{equation}
Thus,
\begin{equation}
\begin{aligned}
&\expect\left[\sum_{i=1}^{\nm}\left((Af^*)_i^{1/2}-(A\wh{f})_i^{1/2}\right)^2\right] \le \min_{(D^\T f)\in\Gamma}\left[2\nm\|A(f^*-f)\|_2^2 + \frac{1}{\inty}\pen(f)\right].
\end{aligned}
\label{eq:upperc}
\end{equation}
Now we can apply Assumption~\ref{as:Arip1} and have
\begin{equation}
\begin{aligned}
\|A (f^*-f)\|_2^2 = \frac{1}{(\resc)^2}\|\wt{A}(f^*-f)\|_2^2 \le \frac{\Ru}{4(a_u-a_{\ell})^2\nm} \|f^*-f\|_2^2, \quad \forall (D^\T f)\in\Gamma .
\end{aligned}
\label{eq:uppere}
\end{equation}
Combining Eq.~\eqref{eq:uppera}, \eqref{eq:upperb}, \eqref{eq:upperc}, and \eqref{eq:uppere}, we get 
\begin{equation}
\begin{aligned}
\expect\left[\|f^*-\wh{f}\|_2^2\right] 	&\le \frac{(\resc)^2}{\Rl}\cdot\frac{4}{\nm}
	\min_{(D^\T f)\in\Gamma}\left[\frac{\Ru}{2(a_u-a_{\ell})^2}\|f^*-f\|_2^2 + \frac{1}{\inty}\pen(f)\right] \\
	&={\min_{(D^\T f)\in\Gamma }\left[\frac{8\Ru}{\Rl}\|f^*-f\|_2^2 
	+ \frac{16(a_u-a_{\ell})^2}{\Rl\inty}\pen(f)\right]}.
\end{aligned}
\label{eq:upper_ori}
\end{equation}
We now bound $\|f^* - f\|_2^2$ and $\pen(f)$ in \eqref{eq:upper_ori} for some good choice of
$f$. In particular, we choose $f_q$ such that $D^Tf_q\in\Theta_K$ is
the quantized version of $\theta^* = D^\T f^*$; \ie $\theta_q$ is the
projection of $\theta^*$ onto $\Theta_K$. For any $f_q\in\cF$, let $\wt{f}_q\in\cC$ be the projection of $f_q$ onto $\cC$, and we have $D^\T\wt{f}_q\in\Gamma$. 
By the Pythagorean identity we have
$\|f^*-\wt{f_q}\|_2^2\le\|f^*-f_q\|_2^2$ because $f^*\in\cC$ and
$\wt{f}_q$ is the projection of $f_q$ onto $\cC$. Furthermore, we have $\pen(\wt{f}_q) = \pen(f_q)$ because we use the same code for $f_q\in \cF$ and its projection $\wt{f}_q \in \cC$. Thus
\begin{equation}
\begin{aligned}
\expect\left[\|f^*-\wh{f}\|_2^2\right]
	\le \frac{8\Ru}{\Rl}\|f^*-{f}_q\|_2^2 
	+ \frac{16(a_u-a_{\ell})^2}{\Rl\inty}\pen({f}_q).
\end{aligned}
\label{eq:tmp1}
\end{equation}
We now bound $\|f^* - f_q\|_2^2$ in \eqref{eq:tmp1}.
Since $\theta_q\in\Theta_K$ is the $(s+1)$-sparse quantized version
of $\theta^*$ so that $(\theta_q)_1 = \theta^*_1$ and $|(\theta_q)_i -
\theta^*_i| \le \frac{L}{K}, i = 2,\ldots,\nc$, we have
\begin{equation}
\begin{aligned}
\|f^* - f_q\|_2^2 = \|\theta^* - \theta_q\|_2^2 \le s\left(\frac{L}{K}\right)^2 = \frac{sL^2}{K^2}.
\end{aligned}
\label{eq:quanerror}
\end{equation}
Using \eqref{eq:defPen}, we have 
\begin{equation}
\begin{aligned}
\pen(f_q) &= 2\log_2s + 2s\log_2\nc + 2s\log_2 K.
\end{aligned}
\notag
\end{equation}
Thus we can bound 
\begin{equation}
\begin{aligned}
\expect\left[\|f^*-\wh{f}\|_2^2\right] \le \frac{8\Ru}{\Rl}\frac{{s}L^2}{K^2} 
+ \frac{32(a_u - a_\ell)^2}{\Rl\inty}[\log_2s + {s}\log_2\nc + {s}\log_2 K],
\end{aligned}
\label{eq:newT1}
\end{equation}
which depends on $K$. The value of $K^*$ defined in \eqref{eq:K}
minimizes this bound, and note that $1<K*\le K\le K^*+1$.
Thus
\begin{equation}
\begin{aligned}
\expect\left[\|f^*-\wh{f}\|_2^2\right] 
&\le \frac{8\Ru}{\Rl}\frac{{sL^2}}{(K^*)^2} 
+ \frac{32(a_u - a_\ell)^2}{\Rl\inty}[\log_2s + {s}\log_2\nc + {s}\log_2 (K^*+1)] \\
& \le \frac{8\Ru}{\Rl}\frac{{sL^2}}{(K^*)^2} 
+ \frac{32(a_u - a_\ell)^2}{\Rl\inty}[\log_2s + {s}\log_2\nc + {s}\log_2 \left((K^*)^2+1\right)].
\end{aligned}
\notag
\end{equation}
Note that $L \le 1$ because $D$ is orthonormal, and using the value for $K^*$ and $L$ yields 
\begin{equation}
\begin{aligned}
\expect\left[\|f^*-\wh{f}\|_2^2\right] 
&\le \frac{16(a_u-a_{\ell})^2s}{\Rl\inty\log2} + \frac{32(a_u-a_{\ell})^2}{\Rl }\frac{1}{\inty}[\log_2s + {s}\log_2\nc] \\
& + \frac{32(a_u-a_{\ell})^2}{\Rl }\frac{s}{\inty}
\left[\log_2\left({\frac{\Ru\log2}{2({a_u-a_\ell})^2}}\inty +1\right) \right].
\end{aligned}
\notag
\end{equation}
Note there exists some constant $c_0>0$ s.t. $\log_2\left({\frac{\Ru\log2}{2({a_u-a_\ell})^2}}\inty +1\right)\le c_0 \log(\inty+1)$, thus the dominating terms in the bound are proportional to
$\frac{s}{\inty}\log_2\nc$ and $\frac{s}{\inty}\log_2(\inty+1)$. Then
there exist some constants $c_1, c_2>0$ such that 
\begin{equation}
\begin{aligned}
\expect[\|f^* - \wh{f}\|_2^2]
	& \le \frac{c_1 {s}\log_2\nc}{\Rl\inty} + \frac{c_2 {s}\log_2(\inty+1)}{\Rl\inty}.
\end{aligned}
\label{eq:up1}
\end{equation}
Instead of using the quantized estimator $\wh{f}$, an alternative
estimator is the mean estimator $\wh{f} =
\frac{1}{\sqrt{\nc}}\bI_{\nm\times1}$. This estimator corresponds to
the case where $K=1$ above.
\begin{equation}
\begin{aligned}
\expect[\|f^* - \wh{f}\|_2^2] = \|f^* - \frac{1}{\sqrt{\nc}}\bI_{\nm\times1}\|_2^2 = \|\btheta^*\|_2^2.
\end{aligned}
\label{eq:up2}
\end{equation}
$\|\btheta^*\|_2^2$ can be upper bounded in one of two ways. Either $\|\btheta^*\|_2^2 <\|f^*\|_2^2\leq \|f^*\|_1^2  = 1$, or
\begin{equation}
\|\btheta^*\|_2^2 \leq s \|\theta^*\|_{\infty}^2  = s \|D^T f^*\|_{\infty}^2 \leq s \max_{j,k}|D_{j,k}|^2 \|f^*\|_1^2 = s \max_{j,k}|D_{j,k}|^2. 
\notag
\end{equation}
Since the minimax mean-squared error takes the minimum over all measurable estimators, we have for some $C_U>0$,
\begin{eqnarray}
\min_{\wh{f}} \max_{f^* \in \FunClass} \expect[\|f^* - \wh{f}\|_2^2]
	& \le C_U \max \big(1, s \max_{j,k}|D_{j,k}|^2,  \frac{s \log_2\nc}{\Rl\inty} + \frac{s \log_2(\inty+1)}{\Rl\inty} \big).
\label{eq:up1}
\end{eqnarray}
\end{proof}

As noted earlier, the lower bound in Theorem~\ref{thm:lower} only requires Assumptions~\ref{as:wtA} and~\ref{as:Arip1}, whereas the upper bound requires the additional Assumption~\ref{as:Arip2}. To explain this, first note that by Lemma~\ref{lem:AfBound}, the intensity of $Af^*$ is both lower and upper bounded. Then with the upper bound of $\|\wt{A}f\|_2^2$, the KL-divergence is upper bounded by the squared $\ell_2$ distance; with the lower bound of $\|\wt{A}f\|_2^2$, the KL-divergence is lower bounded by the squared $\ell_2$ distance. The lower bound only requires that the KL-divergence between parameters in the space is upper bounded by the squared $\ell_2$ distance whereas the upper bound requires that the KL-divergence be both upper and lower bounded in terms of $\ell_2$ distance.

\section{Conclusion}
\label{sec:Conclusion}
\label{sec:Conclusion}
In this paper, we provide sharp reconstruction performance bounds for
photon-limited imaging systems with real-world constraints. These
critical physical constraints, referred to as the positivity and
flux-preserving constraints, arise naturally from the fact that we
cannot subtract light in an optical system, and that the total number
of available photons (which is proportional to the total observation
time or total source intensity) is fixed regardless of the number of sensors in a system. These constraints, which are often not considered in other literature, play an essential role in our performance guarantees.

Unlike most compressed sensing results, our mean-squared error bounds do not decrease as $\nm$, the number of observations, grows to infinity. This independence of $\nm$ derives from the flux-preserving constraint. In most compressed sensing results, the signal-to-noise ratio is controlled by the number of observations, while here, as we work in a photon-limited environment,  we cannot afford an infinite photon-count or observation time. In this photon-limited framework, the signal-to-noise ratio is controlled by the total intensity, $\inty$. This allows us to separate the effect of having a higher signal-to-noise ratio (larger $\inty$) and the effect of having a better variety of observations (larger $\nm$).  As a result, we see that we cannot drive the reconstruction error to 0 simply by having more measurements. Instead, we need infinite signal-to-noise ratio ($\inty\rightarrow \infty$) to achieve perfect reconstruction. However, it is still critical to have $\nm$ sufficiently large in order to have all the assumptions on the sensing matrix satisfied.

Another interesting observation is that the lower and upper bounds
derived in this paper both exhibit different behaviors in
 the low- and 
high-intensity regimes.  Such behavior change is also
demonstrated in the simulations, where we see ``elbows'' in the plot
of MSE vs.\ $\inty$. The low-intensity region, in particular,
corresponds to a range of $\inty$ where reliable reconstruction is
hard to achieve. 

In conventional CS settings, CS does not lead to significantly higher
MSEs than directly sensing nonzero coefficients (if their locations
were known). In contrast, because of the unusual role of noise in our
setting, directly sensing nonzero coefficients can lead to dramatic
reductions in MSE. This is the reason that optical systems which
measure low-resolution images directly rather than collect compressive
measurements can perform so much better in practice, particularly in
low-intensity regimes, as detailed in Section~\ref{sec:dnsamp}. For
example, most night vision cameras have a focal plane array with $n$
elements, each of which directly measures a distinct region of the
field of view. Conventional CS suggests we could change the optical
design so that the same $n$-element focal plane array could be used to
generate an estimate of the scene with $p > n$ pixels. However, our
results suggest that this approach will not be successful unless $T$
(the observation/stare time or scene brightness) is quite
large.

Finally, we note that our minimax upper bounds do not correspond to a feasible
estimator. Developing upper bounds for estimators which are
implementable in polynomial time (\eg corresponding to $\ell_1$
regularization as in the LASSO) is an important avenue for future
work. Because of the close correspondence between our theoretical
upper bounds and our simulations (which used an $\ell_1$
regularization algorithm), we anticipate that the bounds on such a
method would not differ significantly from the bounds presented in
this paper.

\section{Acknowledgements}
We gratefully acknowledge the support of the awards AFOSR
FA9550-11-1-0028 and NSF CCF-06-43947. 
GR was partially supported from August 2012 to July 2013 by
the National Science Foundation under Grant DMS-1127914 to the
Statistical and Applied Mathematical Sciences Institute.
We also thank Dr. Zachary Harmany for the fruitful discussions.

\appendix

\section{Proofs}
\subsection{Proof of Lemma~\ref{lem:flux}}
\label{sec:proofflux}
\begin{proof}
By $\wt{A}_{i,j}\in[a_\ell, a_u]$, we have
\begin{equation}
\begin{aligned}
A_{i, j} = \frac{\wt{A}_{i,j}+\shft\bI_{\nm\times\nc}}{\resc}
\in \left[\eps,\frac{1}{{\nm}}\right].
\end{aligned}
\label{eq:rangeofA}
\end{equation}
Thus we have,
\begin{equation}
\begin{aligned}
A_{i,j} &\ge 0, &\forall i,j,  \\
\sum_{i=1}^n A_{i,j} &\le \sum_{i=1}^n \frac{1}{n} = 1, &\forall j.
\end{aligned}
\notag
\end{equation}

\end{proof}
The value of shift and rescaling parameters in
Eq.~\eqref{eq:Aconstruct} are chosen to ensure
$A_{i,j}\in[\frac{1}{2\nm},\frac{1}{\nm}]$. The upper bound of $1/\nm$
is needed to satisfy the flux-preserving constraint, while the lower
bound of $1/(2\nm)$ is chosen to ensure that that we have a reasonable
lower bound on $Af^*$ so that we can lower bound the KL divergence
between $y$ and $TAf^*$. It is possible to change the shift
and rescaling parameters to have an arbitrary lower bound
$\epsilon<1/\nm$ on $A_{i,j}$. However, $\epsilon$ being either too
small ($\approx 0$) or too big ($\approx 1/\nm$) will result in
suboptimal bounds. The best choice of $\epsilon$ is $c/\nm$ for some
$0<c<1$, where we use $c=1/2$ in Assumption~\ref{as:wtA}. 

\subsection{Proof of Lemma~\ref{lem:tightLB}}
\label{sec:prooftightLB}
The proof technique is similar to the proof in Appendix~\ref{sec:proofLower}. We omit some repeated steps for conciseness.
\begin{proof}
First define the sets
\begin{equation}
\begin{aligned}
\cG &= \{f\in\{0,1\}^{\nc}: \|f\|_0=1\},\\
\cH &= \{\theta\in\reals^{\nc}: \theta = D^\T f, f\in \cG\}.
\end{aligned}
\notag
\end{equation}
Note $\cG\subset\FunClass$, and for this Haar Wavelet basis, all $\theta\in\cH$ satisfies $\|\theta\|_0=1+\log_2\nc$. Thus $\theta$ is $(s+1)$-sparse for all $s\ge\log_2\nc$. Also we have
\begin{equation}
\|\si-\sj\|_2^2 = \|f^i-f^j\|_2^2 = 2, \forall i,j=1,\ldots,M.
\notag
\end{equation}
Thus $\cH$ is a $\sqrt{2}$-packing set of $\FunClass$ in the $\ell_2$ norm. The cardinality of $\cH$ is $M \deq |\cH| = |\cG| = \nc$. Let the elements of $\cH$ be denoted $\{\theta^1,\ldots,\theta^M\}$ and $f^i\deq D\theta^i, i=1,\ldots,M$. Then we can define a random vector $\wt{\Theta}\in\reals^\nc$ that is drawn from a uniform distribution over the set $\cH$, and form a multi-way hypothesis testing problem where $\wt{\theta}$ is the testing result that takes value in $\cH$. Then we can bound the minimax estimation error using the same proof technique in Appendix~\ref{sec:proofLower} and get
\begin{equation}
\begin{aligned}
\min_{f}\max_{f^*\in\FunClass} \expect\|f-f^*\|_2^2 
&\ge \frac{1}{2}\min_{\wt{\theta}}\prob[\wt{\theta}\neq\wt{\Theta}]\\
&\ge \frac{1}{2}\left(1-\frac{I(y;\wt{\Theta}) + \log{2}}{\log{M}}\right),
\end{aligned}
\notag
\end{equation}
where
\begin{equation}
\begin{aligned}
I(y;\wt{\Theta})
&\le \frac{1}{\binom{M}{2}}\sum_{\substack{i,j=1,\ldots,M \\ i\neq j}}\KL(p(y|\inty\Phi\si)\|p(y|\inty\Phi\sj))\\
&\le \frac{\Ru\inty}{2(a_u-a_\ell)^2} \|(\si-\sj)\|_2^2\\
&= \frac{\Ru\inty}{(a_u-a_\ell)^2}.
\end{aligned}
\notag
\end{equation}
By $\log M = \log \nc \ge 2\log2$ and $T\le \frac{(a_u-a_{\ell})^2\log\nc}{4\Ru}$, we have 
\begin{equation}
1-\frac{I(y;\wt{\Theta}) + \log{2}}{\log{M}}\ge \frac{1}{4}.
\notag
\end{equation}
Thus there exists an absolute constant $C_L''\ge 1/8$ s.t. 
\begin{equation}
\min_{f}\max_{f^*\in\FunClass} \expect\|f-f^*\|_2^2 \ge C_L''.
\notag
\end{equation}

\end{proof}

\subsection{Proof of Lemma~\ref{prop:dsUpper}}
\label{sec:proofDS}
\begin{proof}
Recall that we know the first coefficient of $\theta^*$. The total mean square error for the downsampling method is then
\begin{equation}
\begin{aligned}
R^{{\rm DS}} 
	&= \expect\left[\|\wh{\theta}^{{\rm DS}}-\theta^*\|_2^2\right] \\
	&= \sum_{i=2}^{\nc/\kappa} \expect\left[(\wh{\theta}^{{\rm DS}}_i-\theta^*_i)^2\right] + \sum_{j=\nc/\kappa+1}^{\nc} \expect\left[(\wh{\theta}^{{\rm DS}}_j-\theta^*_j)^2\right].
\end{aligned}
\label{eq:ds_decomp}
\end{equation}
In Eq~\eqref{eq:ds_decomp}, $i=2,\ldots,\nc/\kappa$ correspond to the coefficients in the coarse scales. For DWT, the downsampling process does not distort these coefficients, and the reconstruction is unbiased, \ie  
\begin{equation}
\expect\left[\wh{\theta}^{{\rm DS}}_i-\theta^*_i\right]=0,\quad i=2,\ldots,\nc/\kappa.
\notag
\end{equation}
Thus by bias-variance decomposition, $\sum_{i=2}^{\nc/\kappa}\expect\left[(\wh{\theta}^{{\rm DS}}_i-\theta^*_i)^2\right]$ only depends on the variance of the estimator $\wh{\theta}^{{\rm DS}}_i$. To compute ${\rm Var}\left[\wh{\theta}^{{\rm DS}}_i\right]$, we first evaluate ${\rm Var}\left[\wh{f}^{{\rm DS}}_k\right]$ for $k=1,\ldots,\nc$, and then use the basic property of the variance to compute the variance of $\wh{\theta}^{{\rm DS}}_i=D^\T \wh{f}^{{\rm DS}}$. Note that $A^\T $ has only one non-zero entry (with value 1) in each row, thus each of the estimator $\wh{f}^{{\rm DS}}_k$ is a scaled Poisson random variable, whose variance is the scaled mean parameter of the Poisson process:
\begin{equation}
{\rm Var}\left[\wh{f}^{{\rm DS}}_k\right]=\frac{1}{\kappa^2\inty^2}(\inty Af^*)_{k'},
\notag
\end{equation}
where $k'=\lfloor k/\kappa \rfloor +1$, and $(Af^*)_{k'} = \sum_{j=1}^\kappa f^*_{(k'-1)\kappa+j}$. Note that $f^*_j\le 2/{\nc}$ because
\begin{equation}
\begin{aligned}
f^*_i	& = (D\theta^*)_i \\
	& = D_{i,1}\theta_1 + \sum_{j=2}^{\nc}D_{i,j}\theta_j \\
	& = \frac{1}{\sqrt{\nc}}\theta_1 + \alpha\langle \bD_i, \bar{\theta} \rangle \\
	& \le \frac{1}{\sqrt{\nc}}\frac{1}{\sqrt{\nc}} + \alpha_k |\langle \bD_i, \bar{\theta}\rangle| \\
	& \le \frac{1}{\nc} + \alpha \lambda \\
	& = \frac{1}{\nc} + \frac{1}{\nc\lambda}\lambda \\
	& = \frac{2}{\nc},
\end{aligned}
\notag
\end{equation}
and we get
\begin{equation}
{\rm Var}\left[\wh{f}^{{\rm DS}}_k\right] \le \frac{1}{\kappa^2\inty^2}\kappa \frac{2\inty}{\nc} = \frac{2}{\inty\nc\kappa}.
\notag
\end{equation}
Now by $\wh{\theta}^{{\rm DS}}_i=D^\T \wh{f}^{{\rm DS}}$, and $D$ orthonormal, we have
\begin{equation}
\begin{aligned}
{\rm Var}\left[\wh{\theta}^{{\rm DS}}_i\right] 
	&= \sum_{k=1}^{\nc} D_{ki}^2{\rm Var}\left[\wh{f}^{{\rm DS}}_k\right] \\
	&\le \frac{2}{\inty\nc\kappa}\sum_{k=1}^{\nc} D_{ki}^2 \\
	&= \frac{2}{\inty\nc\kappa}.
\end{aligned}
\notag
\end{equation}
Note we have $s'$ coefficients $\theta_i$ that are in the coarse scale, we can conclude by
\begin{equation}
\begin{aligned}
\sum_{i=2}^{\nc/\kappa} \expect\left[(\wh{\theta}^{{\rm DS}}_i-\theta^*_i)^2\right]
	&\le \frac{2s'}{\inty\nc\kappa}.
\end{aligned}
\notag
\end{equation}

On the other hand, $j=\nc/\kappa+1,\ldots,\nc$ correspond to the coefficients in the fine scales. For the DWT basis, the downsampling process completely erases these information because all details in the fine scales become unidentifiable after the downsampling process. As a result, $\wh{\theta}^{{\rm DS}}_j=0$, and $\expect\left[(\wh{\theta}^{{\rm DS}}_j-\theta^*_j)^2\right] = (\theta^*_j)^2 = \alpha_{{s}}^2 = \frac{1}{\nc^2\lambda^2}$.

Add the coarse scale and fine scale coefficients together, we can get
\begin{equation}
\begin{aligned}
R^{{\rm DS}} 
	&= \expect\left[\|\wh{\theta}^{{\rm DS}}-\theta^*\|_2^2\right] \\
	&= \frac{2s'}{\inty \nc \kappa} + \frac{{s}-s'}{\lambda^2\nc^2},
\end{aligned}
\notag
\end{equation}
which completes the proof.
\end{proof}

\section{Others}
\label{sec:calculation}

\subsection{Generation of signals with different sparse supports}
\label{sec:SparseSupportDetail}
$\theta^*_1$ is the triangular signal, where
\begin{equation}
\begin{aligned}
\left(\theta^*_1\right)_j = \begin{cases}
\frac{1}{\sqrt{\nc}}, & j= 1\\
\frac{(s-j+2)}{s\sqrt{\nc}}, & j = 2,\ldots,s+1\\
0, &\mbox{otherwise}.
\end{cases}
\end{aligned}
\notag
\end{equation}
It corresponds to the DCT coefficients of a (rescaled) $\mbox{sinc}^2$ function. $\theta^*_2$ is the signal that is $(s+1)$-sparse with uniform-value non-zero coefficients and the location of the non-zeros being random. In particular, $\theta^*_2$ belongs to the packing set $\cH_{\alpha_s}$ with DCT basis.
$\theta^*_3$ is:
\begin{equation}
\begin{aligned}
\left(\theta^*_3\right)_j = \begin{cases}
\frac{1}{\sqrt{\nc}}, & j=1\\
\frac{2^{k/2}}{\sqrt{\nc}}, & j = 2^k+1, k = 0,\ldots,\min(s,\log_2\nc)-1\\
0, &\mbox{otherwise}.
\end{cases}
\end{aligned}
\notag
\end{equation}
When $s\ge \log_2\nc$, $\theta^*_3$ corresponds to the DWT coefficients of signal $[1,0,\ldots,0]^\T$.
$\theta^*_4$ is the signal in the packing set $\cH_{\alpha_s}$ with DWT basis, where each signal has $s$ randomly located uniform-value coefficients.

\subsection{Calculation of $\lambda_k$}
\label{sec:calcLambda}

\subsubsection{Discrete cosine transform (DCT)}

For a DCT transform matrix $D^{\rm DCT}$,
\begin{equation}
\begin{aligned}
D^{\rm DCT}_{i,j}=\sqrt{\frac{2}{\nc}}\cos{\left(\frac{(2i-1)(j-1)\pi}{2\nc}\right)}, i,j=1,\ldots,\nc.
\end{aligned}
\label{eq:dctval}
\end{equation}
For $1\le k\le s$, we can upper bound $\lambda_{k}(\bD)$ and $\lambda'(\bD)$ using
\begin{equation}
\begin{aligned}
\lambda_{k}\left(\bD^{\rm DCT}\right) 		
				\le& \sqrt{\frac{2}{\nc}}\sum_{j=2}^{k+1}1 
				= \frac{\sqrt{2}k}{\sqrt{\nc}}.
\end{aligned}
\notag
\end{equation}
The upper bound for $\lambda_{k}$ is tight when $\nc$ is large. To see this, fix $i=1$, so we get a lower bound for the two quantities:
\begin{equation}
\begin{aligned}
\lambda_{k}\left(\bD^{\rm DCT}\right) 			
				&\ge \sqrt{\frac{2}{\nc}}\sum_{j=2}^{k+1}\left|\cos\left(\frac{(j-1)\pi}{2\nc}\right)\right| 
				\approx
				\sqrt{\frac{2}{\nc}}\sum_{j=2}^{k+1} \cos 0
				= \frac{\sqrt{2}k}{\sqrt{\nc}},
\end{aligned}
\notag
\end{equation}
when $\nc\gg k$.

\subsubsection{Discrete Hadamard transform (DHT)}

The discrete Hadamard transform matrix $D^{\rm DHT}_{(m)}\in\reals^{2^m\times 2^m}$ has the form
\begin{equation}
\begin{aligned}
\left(D^{\rm DHT}_{(m)}\right)_{i,j}=\frac{1}{2^{\frac{m}{2}}}(-1)^{i\cdot j}.
\end{aligned}
\label{eq:dhtval}
\end{equation}
Thus, because $\nc = 2^m$, we have
\begin{equation}
\begin{aligned}
\lambda_{k}\left(\bD^{\rm DHT}_{(m)}\right) 	
			= \frac{1}{2^{\frac{m}{2}}}\sum_{i=2}^{k+1}1 = \frac{k}{2^{\frac{m}{2}}}
			= \frac{k}{\sqrt{\nc}}.
\end{aligned}
\notag
\end{equation}

\subsubsection{Discrete Haar wavelet transform (DWT)}

For a Haar wavelet transform matrix of dimension $2^m$ by $2^m$, the non-zeros have magnitudes in the set $\{2^{-m/2}, 2^{-(m-1)/2}, \ldots, 2^{-1/2}\}$. Thus $\lambda_k$ will be a sum of a geometric sequence (let $m' = \min{(k,m-1)}$):

\begin{equation}
\begin{aligned}
\lambda_{k}(\bD^{\rm DWT}_{(m)}) &= \frac{1}{2^{\frac{m}{2}}}\sum_{i=1}^{m'} 2^{\frac{m-i}{2}} = \frac{1}{2^{\frac{m}{2}}}\left(\frac{2^{\frac{m-1}{2}}(1-(\frac{1}{\sqrt{2}})^{m'})}{1-\frac{1}{\sqrt{2}}}\right)= \frac{1-2^{-\frac{m'}{2}}}{\sqrt{2}-1}  \rightarrow \frac{1}{\sqrt{2}-1},
\end{aligned}
\notag
\end{equation}
as $m'$ goes to infinity.
\subsection{Comparison of compressed sensing to the downsampling method}
\label{sec:dscs_detail}
\subsubsection{Setup}
Assume $f^*$ is sparse under the DWT basis, and $\theta^*$ has exactly $s+1$ non-zero elements. Again, $s$ is the number of non-zeros in the non-DC coefficients. Note that in order to make the downsampling scheme to work well, we need the signal to bear a certain amount of smoothness. We separate the non-DC coefficients into coarse scales ($\theta_i^*, i=2,\ldots,p/\kappa$, low-frequency coefficients) and fine scales ($\theta_i^*, i=p/\kappa+1,\ldots,p$, high-frequency coefficients). The downsampling process will discard all of the information in the fine scales, but directly sense the coarse-scale coefficients. 
Let $s'$ be the number of non-zeros in coarse scales , $\lambda$ be our sparse-localization constant ($\approx\frac{1}{\sqrt{2}-1}$ and does not depend on ${s}$ for DWT, so we omit the subscript here). We further assume that all non-zero non-DC coefficients have the same magnitude so that the signal energy in the coarse scales is non-trivial and can be reflected by $s'$. Specifically, we use signals from the packing set $\cH_{{s},\alpha_{{s}}}$ where $\alpha_{{s}} = \frac{1}{\nc\lambda}$.

Let $A^{{\rm DS}}$ and $A^{{\rm CS}}$ be the sensing matrices using downsampling and compressed sensing methods defined in Eq.~\eqref{eq:dsmatrix} and Eq.~\eqref{eq:Aconstruct}~\eqref{eq:A0construct2}, respectively; $y^{{\rm DS}}\sim {\rm Poisson}(A^{{\rm DS}} f^*)$ and $y^{{\rm CS}}\sim {\rm Poisson}(A^{{\rm CS}} f^*)$ be the corresponding observations. Let $\wh{\theta}^{{\rm DS}}$ and $\wh{\theta}^{{\rm CS}}$ be the reconstructions using downsampling and compressed sensing methods, respectively. For compressed sensing, we use the penalized maximum likelihood estimator which has been introduced in Section~\ref{sec:UpperBound}. 

Similar to our main results, we evaluate the MSE which is defined as
\begin{equation}
\begin{aligned}
R^{{\rm DS}} &\deq \expect\left[\|\wh{\theta}^{{\rm DS}}-\theta^*\|_2^2\right], \\
R^{{\rm CS}} &\deq \expect\left[\|\wh{\theta}^{{\rm CS}}-\theta^*\|_2^2\right].
\end{aligned}
\notag
\end{equation}

\subsubsection{Upper bound of downsampling}
In the following lemma we provide an upper bound of the downsampling method for the set of signals we use in the simulations.
\begin{lemma}[Upper bound of the downsampling method]
\label{prop:dsUpper}
Let $D$ be the DWT basis. For all $\theta^*\in\cH_{{s},\alpha_{{s}}}$ with $\alpha_{{s}} = \frac{1}{\nc\lambda}$, if there are $s'\le s$ non-zeros in the coarse scales, and we use the sensing/reconstruction method defined in Eq.~\eqref{eq:dsmatrix} and~\eqref{eq:dsrecover}, then
\begin{equation}
\begin{aligned}
R^{{\rm DS}} \le \frac{2s'}{\inty \nc \kappa} + \frac{{s}-s'}{\lambda^2\nc^2}.
\end{aligned}
\notag
\end{equation}
\end{lemma}
The proof is provided in Appendix~\ref{sec:proofDS}.

For the same set of signal $\cH_{{s},\alpha_{{s}}}$, assume we are in the low-intensity region, where the lower bound gives
\begin{equation}
\begin{aligned}
R^{{\rm CS}} \ge \frac{{s}}{\lambda^2\nc^2}.
\end{aligned}
\label{eq:cs_const}
\end{equation}
Thus the difference is at least 
\begin{equation}
\begin{aligned}
R^{{\rm CS}} - R^{{\rm DS}} \ge \frac{s'}{\lambda^2\nc^2} -
\frac{2s'}{\inty \nc \kappa},
\end{aligned}
\notag
\end{equation}
which is larger than zero (downsampling performs better) when $\inty/\nc > c/\kappa$ for some constant $c$. Also, note that this difference is proportional to $s'$, which means when the signal energy in the coarse scales is larger, the performance difference is also larger.

\bibliography{Biblio_LL}
\end{document}